\documentclass[12pt]{article}
\usepackage{amsfonts}
\usepackage{amssymb}
\newtheorem{theorem}{Theorem}

\newtheorem{proposition}[theorem]{Proposition}
\newtheorem{lemma}[theorem]{Lemma}
\newtheorem{corollary}[theorem]{Corollary}
\newtheorem{remark}[theorem]{Remark}

\newcommand{\R}{\mathbb{R}}
\newcommand{\Q}{\mathbb{Q}}
\newcommand{\Sf}{\mathbb{S}}

\newcommand{\Les}{\mathbb{L}}

\newcommand{\Hy}{\mathbb{H}}

\newcommand{\Oes}{\mathbb{O}}

\newcommand{\spa}{\mbox{span}}
\newcommand{\hess}{\mbox{Hess\,}}

\newcommand{\grad}{\mbox{grad}}

\newcommand{\po}{{\hspace*{-1ex}}{\bf .  }}
\newcommand{\ii}{isometric immersion }

\newcommand{\rtf}{Ribaucour transform }

\def\Fes{{\cal F}}

\def\<{\langle}

\def\>{\rangle}
\def\a{\alpha}
\def\va{\varphi}

\def\d{\partial}
\def\bea{\begin{eqnarray*} }
\def\eea{\end{eqnarray*} }
\def\be{\begin{equation} }
\def\ee{\end{equation} }

\def\proof{\noindent{\it Proof: }}
\def\qed{\ifhmode\unskip\nobreak\fi\ifmmode\ifinner
\else\hskip5 pt \fi\fi\hbox{\hskip5 pt \vrule width4 pt
height6 pt  depth1.5 pt \hskip 1pt }}

\begin{document}

\title{Blaschke's problem for hypersurfaces}
\author{ Marcos Dajczer and Ruy Tojeiro}
\date{}
\maketitle

\begin{abstract}
We solve Blaschke's problem for hypersurfaces of dimension $n\geq 3$.
Namely, we determine all pairs of Euclidean hypersurfaces $f,\tilde{f}\colon\,M^n\to\R^{n+1}$
 that induce conformal metrics on $M^n$ and envelope a common sphere congruence in $\R^{n+1}$.
\end{abstract}

\section{Introduction}
A fundamental problem in surface theory is to investigate which data are
sufficient to determine a surface in space. For instance, a generic
immersion $f\colon\,M^2\to\R^3$ into Euclidean three-space is determined,
up to a rigid motion, by its induced metric and mean curvature function.
\emph{Bonnet's problem} is to classify all exceptional immersions.
Locally, this was accomplished by Bonnet \cite{bon}, Cartan \cite{ca1} and
Chern \cite{ch}. They split into three distinct classes, namely, constant
mean curvature surfaces, nonconstant mean curvature \emph{Bonnet surfaces}
admitting a one-parameter family of isometric deformations preserving the
mean curvature function, and surfaces that admit exactly one such
deformation giving rise to a so-called \emph{Bonnet pair}. From a global
point of view the problem has been recently taken up by several authors
(see \cite{kpp} and \cite{bob}).

A generic surface $f\colon\,M^2\to\R^3$ is also determined, up to homothety and translation,  by its conformal structure and its Gauss map.
Classifying the exceptions is known as \emph{Christoffel's problem}. All
local solutions were determined by Christoffel himself \cite{chris}.
Besides minimal surfaces, the remaining nontrivial solutions are
\emph{isothermic  surfaces}, which are characterized by the fact that they
admit local conformal parameterizations by curvature lines on the open
subset of nonumbilic points.

Prescribing the Gauss map of a surface $f\colon\,M^2\to\R^3$ can be
thought of as giving a plane congruence (i.e., a two-parameter family of
two-dimensional affine subspaces of $\R^3$) to be enveloped by $f$.
Christoffel's problem can thus be rephrased as finding which surfaces are
not determined by their conformal structure and a prescribed plane
congruence which they are to envelop.

A similar problem in the realm of M\"obius geometry was studied by
Blaschke \cite{bl} and is now  known as Blaschke's problem. It consists in
finding the surfaces that are not determined, up to M\"obius
transformations, by their conformal structure and a given sphere
congruence (i.e, a two-parameter family of spheres) enveloped by them.
Isothermic surfaces show up again as one of the two nontrivial classes of
exceptional cases. An apparently unrelated class appears as the other:
\emph{Willmore surfaces}, which are best known in connection to the
celebrated Willmore conjecture. Willmore surfaces always arise in pairs of
\emph{dual} surfaces, as conformal envelopes of their common \emph{central
sphere congruence}, whose elements have the same mean curvature as that of
the enveloping surfaces at the corresponding points of tangency.

Unlike the case of Willmore dual surfaces, for any  isothermic surfaces
$f,\tilde{f}\colon\,M^2\to\R^{3}$ that arise as exceptional surfaces for
Blaschke's problem the  curvature lines of $f$ and $\tilde{f}$ coincide,
in which case the sphere congruence is said to be \emph{Ribaucour}. Each
element of such a pair is said to be a \emph{Darboux transform} of the
other.

Blaschke's problem was recently studied in \cite{ma} for surfaces of
arbitrary codimension. On the other hand, the investigation of the
analogous to Christoffel's problem for higher dimensional hypersurfaces
$f\colon\,M^n\to\R^{n+1}$, namely, to determine all hypersurfaces that
admit a conformal deformation preserving the Gauss map, was carried out in
\cite{dv}. The isometric version of the problem had been previously solved
in arbitrary codimension in \cite{dg1}.

In this article we solve Blaschke's problem for hypersurfaces: which pairs
of hypersurfaces $f, \tilde{f}\colon\,M^n\to \R^{n+1}$ envelop a common
regular sphere congruence and induce conformal metrics on $M^n$? (see the
beginning of Section $2$ for the meaning of the regularity assumption).
Since pairs of hypersurfaces that differ by an inversion always satisfy
both conditions, they can be regarded as trivial solutions. Thus we look
for nontrivial ones, that is, pairs of hypersurfaces that  do not
differ by a M\"obius transformation of $\R^{n+1}$.

The problem of determining conformal envelopes of \emph{Ribaucour} sphere
congruences was recently treated in arbitrary dimension and codimension in
\cite{to1}, making use of the extension of the Ribaucour transformation
developed in \cite{dt1} and \cite{dt2} to that general setting. They were
named \emph{Darboux transforms} one of each other,  following the standard
terminology of the surface case. However,  the definition in \cite{to1}
does not exclude the possibility of Darboux pairs that differ by a
composition of a rigid motion and an inversion. Thus, here we rule out
from the classification in \cite{to1} the isometric immersions that only
admit such trivial Darboux transforms. Unfortunately, no interesting
higher dimensional analogues of isothermic surfaces arise: in the
hypersurface case, they reduce, up to M\"obius transformations of
Euclidean space, to cylinders over plane curves, cylinders over surfaces
that are cones over spherical curves and rotation hypersurfaces over plane
curves (after excluding the ones that only admit trivial
Darboux transforms). Our main result is that there are no other solutions
of Blaschke's problem for hypersurfaces.

\begin{theorem}\po\label{thm:blaschke}
Let $f, \tilde{f}\colon\,M^n\to \R^{n+1}$, $n\geq 3$, be a nontrivial
solution of Blaschke's problem. Then $f(M)$ and $\tilde{f}(M)$ are, up to
a M\"obius transformation of $\R^{n+1}$, open subsets of one of the
following:

\begin{itemize}
\item[$(i)$] A cylinder over a plane curve.
\item[$(ii)$] A cylinder $C(\gamma)\times\R^{n-2}$, where $C(\gamma)$
denotes the cone over a curve $\gamma$ in $\mathbb{S}^2\subset\R^3$.
\item[$(iii)$] A rotation hypersurface over a plane curve.
\end{itemize}

Conversely, for any hypersurface $f\colon\,M^n\to \R^{n+1}$ that differs
by a M\"obius transformation of $ \R^{n+1}$ from a hypersurface as in
either of the preceding cases there exists $\tilde{f}\colon\,M^n\to
\R^{n+1}$ of the same type as $f$  such that $(f,\tilde{f})$ is a
nontrivial solution of Blaschke's problem. Moreover, $\tilde{f}$ is a
Darboux transform of $f$.
\end{theorem}

 To prove Theorem
\ref{thm:blaschke}, we show that for a pair of hypersurfaces
$(f,\tilde{f})$, that is a solution of Blaschke's problem, the shape
operators are always simultaneously diagonalizable. This reduces the
problem to the  previously discussed case of Ribaucour sphere congruences.
Our approach is as follows. We are first naturally led to study pairs of
conformal hypersurfaces $f, \tilde{f}\colon\,M^n\to \R^{n+1}$, $n\geq 3$,
that satisfy a weaker condition than that of enveloping a common sphere
congruence. In order to describe it, we use  that a sphere congruence in
$\R^{n+1}$ can be regarded as a map $s\colon\,M^n\to \Sf_1^{n+2}$ into the
Lorentzian hypersphere with constant sectional curvature one of Lorentz
space $\mathbb{L}^{n+3}$ (see the beginning of Section $2$ for details).
We study pairs of conformal hypersurfaces $f,\tilde{f}\colon\,M^n\to
\R^{n+1}$, $n\geq 3$, that envelop (possibly different) sphere congruences
$s, \tilde{s}\colon\,M^n\to \Sf_1^{n+2}$ with the same radius function and
which induce the same metric on $M^n$. By the radius function of a sphere
congruence $s\colon\,M^n\to \Sf_1^{n+2}$ we mean the function that assigns
to each point of $M^n$ the (Euclidean) radius of the sphere $s(p)$. We
point out that this condition is not invariant under M\"obius
transformations of Euclidean space. In this way,  we are able to restrict
the candidates of solutions of Blaschke's problem, in the case in which
principal directions are not preserved, to pairs of surface-like
hypersurfaces over surfaces that are solutions of Bonnet's problem in
three-dimensional space forms. We say that a hypersurface
$f\colon\,M^n\to\R^{n+1}$ is \emph{surface-like} if $f(M)$ is 
the image by a  M\"obius transformation of $\R^{n+1}$ of an open subset of one of the following:
\begin{itemize}
\item[$(i)$] a cylinder $M^2\times\R^{n-1}$ over $M^2\subset\R^3$;
\item[$(ii)$] a cylinder $CM^2\times\R^{n-2}$, where $CM^2\subset\R^4$ denotes the cone over a surface
$M^2\subset \Sf^3$;
\item[$(iii)$] a rotation hypersurface over $M^2\subset\R^3_+$.
\end{itemize}
The proof is then completed by showing that none of the possible
candidates is in fact a solution of Blaschke's problem.

\section[Conformally deformable hypersurfaces]{Conformally deformable hypersurfaces}

 Two  hypersurfaces $f\colon M^{n}\to\R^{n+1}$
and $\tilde{f}\colon M^{n}\to\R^{n+1}$ in Euclidean space are said to be
\emph{conformally congruent} if they differ by a conformal transformation
of $\R^{n+1}$. By Liouville's theorem, any such transformation 
is a composition $T=L\circ {\cal I}$ of a similarity $L$ and 
an inversion
${\cal I}$ with respect to a hypersphere of $\R^{n+1}$. Recall that the
inversion $\mathcal{I}\colon\,\R^N\setminus \{p_0\}\to \R^N\setminus
\{p_0\}$ with respect to a hypersphere with radius $r$ centered at $p_0$
is given by
$$
\mathcal{I}(p)=p_0+\frac{r^2}{\|p-p_0\|^2}(p-p_0).
$$
If  $f\colon M^{n}\!\to\!\R^{n+1}$ is a hypersurface and  $N$ is a unit
normal vector field to $f$, then it is easily seen that
$\tilde{N}=r^{-2}\|f-p_0\|^2\mathcal{I}_*N$ defines a unit normal vector
field to $\tilde{f}=\mathcal{I}\circ f$. Moreover, the shape operators
$A_N$ and $\tilde{A}_{\tilde{N}}$ of $f$ and $\tilde{f}$ with respect to
$N$ and $\tilde{N}$, respectively, are  related by \be\label{shape}
r^2\tilde{A}_{\tilde{N}}=\|f-p_0\|^2A_N+2 \<f-p_0,N\>I, \ee where $I$
stands for the identity endomorphism of $TM$. Recall that
$A_NX=-\bar{\nabla}_XN$ for any $X\in TM$, where $\bar{\nabla}$ denotes
the derivative of $\R^{n+1}$.  In particular, $f$ and $\tilde{f}$ have
common principal directions and the corresponding principal curvatures are
related by \be\label{pcurv} r^2\tilde{\lambda}_i=\lambda_i\|f-p_0\|^2
+2\<f-p_0,N\>. \ee

A hypersurface  $f\colon M^{n}\to\R^{n+1}$ is  said to be
\emph{conformally rigid} if  any other conformal immersion
$\tilde{f}\colon\,M^n\to\R^{n+1}$ is conformally congruent to $f$. The
following criterion for conformal rigidity is due to Cartan \cite{ca1}.

\begin{theorem}\po\label{cartan}  A hypersurface
$f\colon\,M^n\to\R^{n+1}$, $n\geq 5,$ is  conformally rigid if all
principal curvatures have multiplicity less than $n-2$ everywhere.
\end{theorem}

 A  conformal immersion $\tilde{f}\colon M^{n}\to\R^{n+1}$ not conformally congruent to $f$ is said to be  a
\emph{conformal deformation} of $f$. It is said to be \emph{nowhere
conformally congruent} to $f$ if it is not conformally congruent to $f$ on
any open subset of $M^n$.

   It is well-known that an
$n$-dimensional Euclidean hypersurface has a principal curvature of
multiplicity at least $n-1$ everywhere if and only if it is conformally flat, hence, highly conformally deformable.  By Theorem~\ref{cartan}, if
an Euclidean hypersurface of dimension \mbox{$n\geq 5$} has principal
curvatures of multiplicity less than $n-1$ everywhere and admits a
conformal nowhere conformally congruent deformation, then it must have a
principal curvature $\lambda$ of constant multiplicity $n-2$ everywhere.
Such a hypersurface was called in \cite{dt3} a \emph{Cartan hypersurface}
if, in addition, $\lambda$ is nowhere zero.

    Cartan hypersurfaces of dimension $n\geq 5$ have been classified in \cite{ca1}.
We refer to  \cite{dt3} for a modern account of that classification as well as for an alternative one. They
can be  separated into four classes, namely,  surface-like, conformally
ruled, the ones having precisely a continuous $1$--parameter family of
deformations and those that admit  only one deformation.

 The approach in \cite{dt3} is based on the structure of
the \emph{splitting tensor} $C$ of the  eigenbundle $\Delta=\ker(A-\lambda
I)$ correspondent to the principal curvature $\lambda$ of multiplicity
$n-2$  of a Cartan hypersurface. It is defined by
$$
\<C_TX,Y\>=\<\nabla_XY,T\>\,\,\,\mbox{for all}\;\;T\in
\Delta\;\;\mbox{and}\;\;X,Y\in \Delta^\perp.
$$
Under a conformal change of metric $\<\,,\,\>^\sim=e^{2\va}\<\,,\,\>$, the
tensor $C$ changes as \be\label{coker}
\tilde{C}_T=C_T-T(\va)I\;\;\mbox{for all}\;\;T\in \Delta. \ee This follows
immediately from the formula \be\label{eq:lcivita}
\tilde{\nabla}_XY=\nabla_XY +X(\varphi)Y+Y(\varphi)X-\<X,Y\>\nabla\varphi,
\ee that relates the Levi-Civita connections $\nabla$ and $\tilde{\nabla}$
of $\<\,,\,\>$ and $\<\,,\,\>^\sim$, respectively. Here $\nabla\varphi $
denotes the gradient with respect to $\<\,,\,\>$.

 A key observation on the splitting tensor  associated to a Cartan hypersurface is the
 following result, which we will also need
here. It slightly improves Lemma $15$ in \cite{dt3}, so we include its
proof.

\begin{lemma}\po\label{le:kerC}
If $C$ is the splitting tensor of $\Delta$, then the dimension of the subspace
${coker}\,C=(\ker C)^\perp$ is at most two at any point of $M^n$.
Moreover, if it is two everywhere then there exists $S\in \mbox{coker}\,C$
such that $C_S=a I$ for some nonzero real number~$a$.
\end{lemma}

\proof It was shown in Lemma $14$ of \cite{dt3} that there exists an
operator $D$ on $\Delta^\perp$  such that $\det D=1$ and $[D,C_T]=0 $ for
all $T\in \Delta$. Thus, the image of $C$ lies in the two-dimensional
subspace $S$ of linear operators on $\Delta^\perp$ that commute with $D$.
This already implies the first assertion. Assuming that the second
assertion does not hold, the subspace spanned by the image of $C$ and the
identity operator would have dimension three and be contained in $S$, a
contradiction.\qed\vspace{1ex}

The simplest structure of the splitting tensor $C$ of a Cartan
hypersurface occurs when there exists a vector field $\delta\in
\Delta^\perp$ such that $C_T=\<\delta,T\>I$ for every $T\in \Delta$. This
is equivalent to requiring $\Delta^\perp$ to be an \emph{umbilical
distribution} with mean curvature vector field $\delta$.
     The following classification  of the corresponding Cartan hypersurfaces was derived in \cite{dt3} as a consequence of the main theorem of  \cite{dft} and plays a key role in this paper.

\begin{theorem}\po\label{umb}  Let $f\colon\,M^n\to\R^{n+1}$, $n\geq 4$, be a
Cartan hypersurface and let $\Delta$  be the eigenbundle correspondent to
its principal curvature of multiplicity $n-2$. If $\Delta^\perp$ is an
umbilical distribution then $f$ is a surface-like hypersurface.
\end{theorem}

The preceding result also holds for $n=3$ if $\lambda$ is assumed to be
constant along $\Delta$, a condition that is always satisfied when the
rank of $\Delta$ is at least two (see \cite{dft}).\vspace{1ex}

 To conclude this section, we point out that conformally
deformable Euclidean hypersurfaces of dimensions $3$ and $4$ have also
been studied by Cartan \cite{ca2},\cite{ca3}, although in these cases a
classification is far from being complete. Even though we do not make use
of Cartan's results for these cases, some of our arguments are implicit in
his work.

\section[ A necessary condition for a solution]{A necessary condition for a solution}

Let $\mathbb{L}^{n+3}$ be the $(n+3)$--dimensional Minkowski space, that
is, $\R^{n+1}$ endowed with a  Lorentz scalar product of signature
$(+,\ldots,+,-)$, and let
$$
\mathbb{V}^{n+2}= \{p\in\mathbb{L}^{n+3}\colon\<p,p\>=0\}
$$
denote the light cone in $\mathbb{L}^{n+3}$. Then
$$
\mathbb{E}^{n+1}=\mathbb{E}^{n+1}_w
=\{p\in\mathbb{V}^{n+2}\colon\<p,w\>=1\}
$$
is a model of $(n+1)$--dimensional Euclidean space for any
$w\in\mathbb{V}^{n+2}$. Namely, choose $p_0\in \mathbb{E}^{n+1}$ and a
linear isometry $D\colon\,\R^{n+1}\to \mbox{span}\{p_0,w\}^\perp\subset
\mathbb{L}^{n+3}$. Then the triple $(p_0,w,D)$ gives rise to an isometry
$\Psi=\Psi_{p_0,w,D}\colon\,\R^{n+1}\to \mathbb{E}^{n+1}\subset
\mathbb{L}^{n+3}$ defined by
$$
\Psi(x)=p_0+Dx-\frac{1}{2}\|x\|^2w.
$$

Hyperspheres can be nicely described in $\mathbb{E}^{n+1}$: given a
hypersphere $S\subset\R^{n+1}$  with (constant) mean curvature $H$ with
respect to a unit normal vector field  $N$ along $S$, then
$v=H\Psi+\Psi_*N\in \mathbb{L}^{n+3}$ is a constant space-like vector.
Moreover, the vector $v$ has unit  length and $\<v,\Psi(q)\>=0$ for all
$q\in S$; thus
$$
\Psi(S)=\mathbb{E}^{n+1}\cap \{v\}^\perp.
$$
Therefore, given a hypersurface $f\colon\,M^n\to \R^{n+1}$, a sphere
congruence enveloped by $f$ with radius function
$R\in\mathcal{C}^\infty(M)$ can be identified with the map
$s\colon\,M^n\to\mathbb{S}_1^{n+2}$ into Lorentzian sphere
$\mathbb{S}_1^{n+2}= \{p\in\mathbb{L}^{n+3}\colon\<p,p\>=1\}$ defined by
\be\label{s} s(q)=\frac{1}{R(q)}\Psi(f(q))+\Psi_*(f(q))N(q). \ee The
sphere congruence is said to be \emph{regular} if the map $s$ is an
immersion.

\begin{proposition}\po\label{prop:indmetrics}
Let $f\colon\,M^n\to\R^{n+1}$, $n\geq 2$, envelop a sphere congruence
$s\colon\,M^n\to\mathbb{S}_1^{n+2}$ with radius function
$R\in\mathcal{C}^\infty(M)$. Then the metrics $\<\,,\,\>$ and
$\<\,,\,\>^*$ induced by $f$ and $s$ are related by \be\label{form}
\<X,Y\>^*=\<(A-\alpha I)X, (A-\alpha I)Y\>, \ee where $\alpha=1/R$ and $A$
is the shape operator of $f$. In particular, the sphere congruence is
regular if and only if $\alpha$ is nowhere a principal curvature of $f$.
\end{proposition}

\proof Differentiating (\ref{s}) we obtain
$$
s_*X=X(\alpha)(\Psi\circ f)+\alpha \Psi_* f_*X-\Psi_*AX-\<N,f_*X\>w.
$$
The conclusion now  follows easily by using that $\<\Psi,\Psi\>=0$, and
hence that $\<\Psi_*Z,\Psi\>=0$ for any $Z\in\R^{n+1}$.\qed\vspace{1ex}

\begin{corollary}\po\label{cor:*}
Let $\tilde{f}, f\colon\,M^n\to\R^{n+1}$, $n\geq 2$, induce conformal metrics
$\<\,,\,\>^\sim=e^{2\va}\<\,,\,\>$ on $M^n$.  Then the following
assertions are equivalent:
\begin{itemize}
\item[$(i)$] $f$ and $\tilde{f}$ envelop sphere
congruences with the same radius function $R$ which induce the same metric
on $M^n$.
\item[$(ii)$] There exists $\alpha\in\mathcal{C}^\infty(M)$ such that $B=A-\alpha I$ and $\tilde{B} =\tilde{A}-\alpha I$ satisfy \be\label{Bs}
\tilde{B}^2=e^{-2\va}B^2.\ee
\end{itemize}
\end{corollary}
\proof By Proposition \ref{prop:indmetrics},  if either $(i)$ or $(ii)$
holds, then so does the other with $\alpha=1/R$. \qed\vspace{1ex}

   Corollary \ref{cor:*} can be extended to pairs of hypersurfaces $f,\tilde{f}\colon\,M^n\to\Q_c^{n+1}$
   in any space form with constant sectional curvature $c$.
If, for simplicity, we take $c=\pm 1$ when $c\neq 0$, then the function
$\alpha$ in part $(ii)$ is related to the radius function $R$ of the
sphere congruence enveloped by $f$ and $\tilde{f}$ by $\alpha=\cot R$ if
$c=1$ and $\alpha=\coth R$ if $c=-1$.

It follows from Corollary \ref{cor:*} that $(ii)$ is a necessary condition
for a pair of conformal hypersurfaces $f,\tilde{f}\colon\,M^n\to\R^{n+1}$
to be a solution of Blaschke's problem, that is, to envelop a common
sphere congruence. This can also be derived directly for hypersurfaces in
$\Q_c^{n+1}$ from the fact that $f$ and $\tilde{f}$ enveloping a common
sphere congruence with radius function $R\in\mathcal{C}^\infty(M)$ is
equivalent to \be\label{env} Cf+SN=C\tilde{f}+S\tilde{N}, \ee where we use
the standard models  $\Sf^{n+1}\subset\R^{n+2}$ and
$\Hy^{n+1}\subset\Les^{n+2}$ of $\Q_c^{n+1}$ when $c=1$ and $c=-1$,
respectively, so that $\<f,f\>=\<\tilde f, \tilde f\>=c$. Moreover,
$$
\left\{ \begin{array}{l} C=\cos R,\;\;\; S=\sin R\;\;\;\;\;\;\mbox{if}\;\;
c=1,
\vspace*{1.5ex}\\
C=1,\;\;\;\;\;\;\;\;\;S=1/R\;\;\;\;\;\;\;\,\mbox{if}\;\; c=0,
\vspace*{1.5ex}\\
C=\cosh R,\;S=\sinh R\;\;\;\;\mbox{if}\;\; c=-1,
\end{array}\right.
$$
whereas $N$ and $\tilde{N}$ are unit vector fields normal to $f$ and
$\tilde{f}$, respectively. Differentiating (\ref{env}) yields
$$
X(R)(Sf+CN)+f_*(CI-SA)X =X(R)(S\tilde f+C\tilde N)+\tilde f_*(CI-S\tilde
A)X.
$$
Setting
$$
\alpha=C/S
$$
this gives \be\label{der} f_*BX-X(R)(f+\a N)=\tilde{f}_*\tilde BX
-X(R)(\tilde f+\a\tilde N), \ee which implies that
$$
\|f_*BX\|=\|\tilde{f}_*\tilde{B}X\|
$$
for any $X\in TM$. It follows that
$$
\<f_*BX,f_*BY\> =\<\tilde{f}_*\tilde{B}X,\tilde{f}_*\tilde{B}Y\>
=e^{2\va}\<{f}_*\tilde{B}X,{f}_*\tilde{B}Y\>
$$
for all $X,Y\in TM$, or equivalently, that
$\tilde{B}^2=e^{-2\va}B^2$.\vspace{1ex}

In the remaining of the present section we study pairs of conformal
hypersurfaces $f, \tilde{f}\colon\,M^n\to\R^{n+1}$ that satisfy condition
$(ii)$ of Corollary \ref{cor:*}.
We point out that the limiting case in which $\alpha$ is identically zero
reduces to the problem recently studied by Vlachos \cite{v} of determining
all pairs of conformal hypersurfaces $f, \tilde{f}\colon\,M^n\to\R^{n+1}$
whose Gauss maps with values in the Grassmannian of $n$-planes in
$\R^{n+1}$ induce the same metric on $M^n$. In particular, this allows us
to adapt to our case some of the arguments used in the proof of the main
result of that paper.

\begin{lemma}\po\label{le:1}
Let $ \tilde{f}, f\colon\,M^n\to\R^{n+1}$, $n\geq 3$, induce the conformal
metrics $\<\,,\,\>^\sim=e^{2\va}\<\,,\,\>$ on $M^n$ and satisfy either one
of the equivalent conditions in Corollary \ref{cor:*}.  Assume that the
shape operators $A$ and $\tilde{A}$ of $f$ and $\tilde{f}$, respectively,
cannot be simultaneously diagonalized at any point of $M^n$. Then there
exist a smooth distribution $\Delta$ of rank $n-2$ such that

\begin{itemize}
\item[(i)] $\Delta$ is the common eigenbundle
$\ker(A-\lambda I)=\ker(\tilde{A}-\tilde{\lambda} I)$ correspondent to
principal curvatures $\lambda$ and $\tilde{\lambda}$ of $f$ and
$\tilde{f}$, respectively,

\item[(ii)] $\mbox{trace } (A|_{\Delta^\perp})=2\alpha =\mbox{trace }(\tilde{A}|_{\Delta^\perp})$,

\item[(iii)] $\ker B|_{\Delta^\perp}=\{0\}=\ker \tilde{B}|_{\Delta^\perp}$,
\end{itemize}

and an orthogonal tensor $T$ on $M^n$ such that

\begin{itemize}
\item[(iv)] $T|_\Delta=\epsilon\,I$ for $\epsilon=\pm 1$,

\item[(v)] $\det T|_{\Delta^\perp}=1$,

\item[(vi)] $\tilde{B}=e^{-\varphi}B\circ T$.
\end{itemize}
\end{lemma}

\proof Let $\lambda_1,\ldots, \lambda_n$ and $\tilde{\lambda}_1, \ldots,
\tilde{\lambda}_n$ be the principal curvatures of $f$ and $\tilde{f}$,
with corresponding principal frames $\{e_1,\ldots, e_n\}$ and
$\{\tilde{e}_1,\ldots, \tilde{e}_n\}$,  respectively, which we assume to
be orthonormal with respect to the metric induced by $f$. Set
$$
\mu_j=\lambda_j-\alpha\;\;\;\;\mbox{and}\;\;\;\;\tilde\mu_j
=\tilde\lambda_j-\alpha.
$$
By condition $(ii)$ of Corollary \ref{cor:*}, after re-enumeration of the
principal vectors, if necessary, we have \be  \label{id}
\tilde{\mu_j}^2=e^{-2\varphi}\mu_j^2 \;\;\;\mbox{for}\;\;\;1\le j\le n.
\ee Write
$$
\tilde{e}_i=\sum_{j=1}^na_{ji}e_j\;\;\; \mbox{for}\;\;\; 1\le i\le n.
$$
Then
$$
\tilde\mu_i^2\sum_j a_{ji}e_j =\tilde{B}^2\tilde{e}_i
=e^{-2\varphi}B^2\tilde{e}_i =e^{-2\varphi}\sum_ja_{ji}\mu_j^2e_j.
$$
Hence,
$$
(\tilde\mu_i^2 - e^{-2\varphi}\mu_j^2)a_{ji}=0 \;\;\;\mbox{for}\;\;\; 1\le
i, j\le n,
$$
and it follows from (\ref{id}) that \be  \label{dif}
(\mu_i^2-\mu_j^2)a_{ji}=0= (\tilde\mu_i^2-\tilde\mu_j^2)a_{ji} \;\;\;
\mbox{for}\;\;\; 1\le i, j\le n. \ee Therefore, if $a_{ji}\neq 0$ then
$\mu_i^2=\mu_j^2$ and $\tilde\mu_i^2=\tilde\mu_j^2$, or equivalently \be
\label{two} (\lambda_i+\lambda_j-2\alpha)(\lambda_i-\lambda_j)=0
\;\;\;\mbox{and}\;\;\; (\tilde\lambda_i+\tilde\lambda_j
-2\alpha)(\tilde\lambda_i-\tilde\lambda_j)=0. \ee

We now assume the existence of a smooth distribution $\Delta$ of rank
$n-2$ satisfying  $(i)$  and prove the remaining assertions. The principal
curvatures can be ordered so that $\Delta=\mbox{span}\, \{e_3,\ldots,
e_n\}= \mbox{span}\,\{\tilde{e}_3,\ldots, \tilde{e}_n\}$,
$\lambda_3=\cdots=\lambda_n:=\lambda$ and $\tilde{\lambda}_3=\cdots=
\tilde{\lambda}_n:=\tilde{\lambda}$. Since $a_{21}\neq 0$ by the
assumption that $A$ and $\tilde{A}$ cannot be simultaneously diagonalized
at any point, then (\ref{two}) is satisfied for $i=1$ and $j=2$. Since
$\lambda_1\neq\lambda_2$ and $\tilde\lambda_1\neq\tilde\lambda_2$ by the
same assumption, we obtain that $(ii)$ holds.

Observe that  $\ker \tilde{B}=\ker B$, in view of  (\ref{Bs}),   thus our
assumption on $A$ and $\tilde{A}$ implies that condition $(iii)$ must
hold.
 Therefore, using (\ref{Bs})
once more, we obtain that all the remaining conditions in the statement
are fulfilled by the tensor $T$ defined by
\begin{itemize}
\item[$(i)$] $ T|_{\Delta^\perp}
=e^\va(B|_{\Delta^\perp})^{-1}\tilde{B}|_{\Delta^\perp};$
\item[$(ii )$]$T|_\Delta=\epsilon I$, where $\epsilon =1$ or $\epsilon =-1$, according as
$\tilde\mu=\tilde{\lambda}-\alpha$ and $\mu=\lambda-\alpha$ are related by
$\tilde\mu=e^{-\varphi}{\mu}$ or $\tilde\mu=-e^{-\varphi}{\mu}$,
respectively.
\end{itemize}

 In order to complete the proof, it suffices to show that
at each $x\in M^n$ there exists a common eigenspace $\Delta(x)$ of $A$ and
$\tilde{A}$ of dimension $n-2$. The assumption on $A$ and $\tilde{A}$
forces $\Delta(x)$ to be maximal with this property, and this implies
smoothness of $\Delta$.  We consider separately the cases $n\geq 5$, $n=3$
and $n=4$.

\vspace{1ex} \noindent \textit{Case $n\ge 5$.} Taking  (\ref{shape}) into
account, it follows from the assumption on $A$ and $\tilde{A}$ that
$f|_U$ and $\tilde{f}|_U$ do not coincide up to a conformal diffeomorphism
of Euclidean space on any open subset $U\subset M^n$. Then, existence of a
common eigenspace $\Delta(x)$ of $A$ and $\tilde{A}$ of dimension $n-2$
follows from Theorem \ref{cartan}. \vspace{2ex}

\noindent \textit{Case $n=3$.} All we have to prove in this case is the
existence of a common principal direction of $A$ and $\tilde{A}$. First
notice that the case in which $\mu_1^2,\mu_2^2,\mu_3^2$ are mutually
distinct is ruled out by (\ref{dif}) and the assumption on $A$ and
$\tilde{A}$. Therefore, we are left with two possibilities: \vspace{1ex}

\noindent $(a)$  $\mu_1^2=\mu_2^2\neq\mu_3^2$, up to a reordering. Then
$a_{13}=a_{23}=0$ by (\ref{dif}), and hence $e_3$ is a common principal
direction of $A$ and $\tilde{A}$.
\medskip

\noindent $(b)$  $\mu_1^2=\mu_2^2=\mu_3^2$. We may assume that
$-\mu_2=\mu_3=\mu_4$ and that $-\tilde\mu_2=\tilde\mu_3=\tilde\mu_4$.
Then, both $A$ and $\tilde A$ have a two dimensional eigenspace and their
intersection yields the desired common principal direction. \vspace{2ex}

\noindent \textit{Case $n=4$.} As before, it follows from (\ref{dif}) and
the assumption on $A$ and $\tilde{A}$ that
$\mu_1^2,\mu_2^2,\mu_3^2,\mu_4^2$ can not be mutually distinct.
\medskip

For the remaining cases we need the following facts. The curvature tensors
$R $ and $\tilde R$ of the metrics induced by $f$ and $\tilde f$ are related by 
\be\label{curvs}
\begin{array}{l}
\tilde R(X,Y)Z=R(X,Y)Z-(Q(Y,Z)+
\<Y,Z\>|\mbox{grad}\varphi|^2)X\vspace{1ex} \\
\hspace*{4ex} + \,(Q(X,Z)+\<X,Z\>|\mbox{grad}\varphi|^2)Y
-\<Y,Z\>Q_0(X)+\<X,Z\>Q_0(Y),
\end{array}
\ee where
$Q_0(X)=\nabla_X\mbox{grad}\varphi-\<\mbox{grad}\varphi,X\>\mbox{grad}\varphi$
and $Q(X,Y)=\<Q_0(X),Y\>$. In particular, if $X,Y,Z$ are orthonormal
vectors then
$$
\tilde R(X,Y)Z=R(X,Y)Z-Q(Y,Z)X+Q(X,Z)Y.
$$
 From the Gauss equation for $f$ it follows that
\be\label{con} 
\<\tilde R(e_r,e_j)e_k,e_r\>=-Q(e_j,e_k)
\;\;\;\mbox{if}\;\;\;r\neq j\neq k\neq r
\ee and \be\label{gauss0} 
\<\tilde R(e_r,e_j)e_k,e_s\>=0\,\,\,\mbox{if
all four indices are distinct}. 
\ee 
In particular, (\ref{con}) implies
that 
\be \label{gauss} \<\tilde R(e_r,e_j)e_k,e_r\>= 
\<\tilde R(e_s,e_j)e_k,e_s\>,\,\,\,\mbox{if 
$\{r,s\}\cap \{j,k\}=\emptyset$ and
$j\neq k$}. \ee

It will be convenient to single out the following consequence of
(\ref{gauss}):\vspace{2ex}\\
{\bf FACT:} \,\, If $e_i=\tilde{e}_i$ and $e_j=\tilde{e}_j$ (up to sign)
for some $1\leq i\neq j\leq 4$  then $\mu_i=\mu_j$
and~$\tilde{\mu}_i=\tilde{\mu}_j$.\vspace{2ex}\\
To prove the Fact, for simplicity of notation we assume 
$(i,j)=(1,2)$. Set
$$
e_3=\cos\theta\tilde e_3 + \sin\theta\tilde e_4,\;\;\;
e_4=-\sin\theta\tilde e_3 + \cos\theta\tilde e_4.
$$
It follows from (\ref{gauss}) for $(r,j,k,s)=(1,3,4,2)$ and the Gauss
equation for $\tilde f$ that
$$
(\tilde\lambda_1-\tilde\lambda_2)(\tilde\lambda_3 -\tilde\lambda_4)
\sin\theta\cos\theta=0
$$
or, equivalently, that
$$
(\tilde\mu_1-\tilde\mu_2)(\tilde\mu_3-\tilde\mu_4) \sin\theta\cos\theta=0.
$$
Since both $\sin\theta\cos\theta=0$ and $\tilde{\mu}_3=\tilde{\mu}_4$ lead
to a contradiction with our assumption on $A$ and $\tilde{A}$, it follows
that $\tilde{\mu}_1=\tilde{\mu}_2$. Reversing the roles of $f$ and
$\tilde{f}$ gives $\mu_1=\mu_2$, and the proof of the Fact is
completed.\medskip

We now proceed with the proof of the existence of a common two-dimensional
eigenspace of $A$ and $\tilde{A}$.
\medskip

\noindent $(a)$ Assume that $\mu_1^2,\mu_2^2,\mu_3^2$ are mutually
distinct and that $\mu_3^2=\mu_4^2$. Then $a_{i1}=0$ if $i\neq 1$ and
$a_{i2}=0$ if $i\neq 2$ by (\ref{dif}). Hence $e_1=\tilde e_1$,
$e_2=\tilde e_2$ (up to sign) and $\mbox{span}\{e_3,e_4\}
=\mbox{span}\{\tilde e_3,\tilde e_4\}$, in contradiction with the
Fact.\medskip

\noindent $(b)$ Assume that $\mu_1^2\neq\mu_2^2=\mu_3^2=\mu_4^2$. Then
$a_{i1}=0$ if $i\neq 1$, and hence $e_1=\tilde e_1$ and
$\mbox{span}\{e_2,e_3,e_4\} =\mbox{span}\{\tilde e_2,\tilde e_3,\tilde
e_4\}$. We may assume $-\mu_2=\mu_3=\mu_4$ and
$-\tilde\mu_2=\tilde\mu_3=\tilde\mu_4$. If $\mbox{span}\{e_3,e_4\}
=\mbox{span}\{\tilde e_3,\tilde e_4\}$ then we are done. Otherwise, we may
assume that $e_4=\tilde e_4$, which gives a contradiction with the Fact.
\medskip

\noindent $(c)$ Assume that $\mu_1^2=\mu_2^2\neq\mu_3^2=\mu_4^2$. Then
$a_{ij}=0$ if $i=1,2$ and $j=3,4$. Hence
$\mbox{span}\{e_1,e_2\}=\mbox{span} \{\tilde e_1,\tilde e_2\}$ and
$\mbox{span}\{e_3,e_4\}=\mbox{span}\{\tilde e_3,\tilde e_4\}$. Since $A$
and $\tilde{A}$ cannot be simultaneously diagonalized we have to consider
only two cases. If $\mu_1=-\mu_2, \mu_3=\mu_4$ and
$\tilde\mu_1=-\tilde\mu_2, \tilde\mu_3=\tilde\mu_4$, then
$\mbox{span}\{e_3,e_4\}=\mbox{span}\{\tilde e_3,\tilde e_4\}$ is the
desired common two dimensional eigenspace of $A$ and $\tilde{A}$.

Assume that $\mu_1=-\mu_2, \mu_3=-\mu_4$ and that
$\tilde\mu_1=-\tilde\mu_2:= \gamma_1\neq 0,
\tilde\mu_3=-\tilde\mu_4:=\gamma_2\neq 0$. Setting
$$
e_1=\cos\phi\tilde e_1 + \sin\phi\tilde e_2,\;\;\; e_2=-\sin\phi\tilde e_1
+ \cos\phi\tilde e_2,
$$
$$
e_3=\cos\theta\tilde e_3 + \sin\theta\tilde e_4,\;\;\;
e_4=-\sin\theta\tilde e_3 + \cos\theta\tilde e_4,
$$
it follows from (\ref{gauss0}) for $(r,j,k,s)=(1,3,4,2)$ and the Gauss
equation for $\tilde f$ that
$$
(\tilde\lambda_1-\tilde\lambda_2)(\tilde\lambda_3 -\tilde\lambda_4)
\sin\theta\cos\theta\sin\phi\cos\phi=0.
$$
Thus $\gamma_1\gamma_2\sin\theta\cos\theta\sin\phi\cos\phi=0$. Hence, we
may assume that $e_1=\tilde e_1$ and $e_2=\tilde e_2$, and the Fact shows
that this case can not occur.
\medskip

\noindent $(d)$ Assume that $\mu_1^2=\mu_2^2=\mu_3^2=\mu_4^2$. If
$-\mu_1=\mu_2=\mu_3=\mu_4$ and $-\tilde\mu_1=\tilde\mu_2=\tilde\mu_3=
\tilde\mu_4$, then both $A$ and $\tilde{A}$ have three-dimensional
eigenspaces, and their intersection gives  a common two-dimensional
eigenspace as desired. If $-\mu_1=\mu_2=\mu_3=\mu_4$ and
$-\tilde\mu_1=-\tilde\mu_2=\tilde\mu_3=\tilde\mu_4$, then we may assume
that $e_2=\tilde e_2$ and $e_3=\tilde e_3$, and we get a contradiction
with the Fact.

To conclude the proof we assume that $ -\mu_1=-\mu_2=\mu_3=\mu_4$ and
$-\tilde\mu_1=-\tilde\mu_2=\tilde\mu_3=\tilde\mu_4$ on an open subset
$U\subset M^4$. Then $f|_U$ and $\tilde{f}|_U$ are Cyclides of Dupin. In
particular, $U$ is conformal to an open subset of a Riemannian product
$\Q_c^2\times \Sf^2$, $c>-1$, with $\Delta_1=\spa\{e_1,e_2\}$ and
$\Delta_2=\spa\{e_3,e_4\}$ as the distributions tangent to the factors
$\Q_c^2$ and $\Sf^2$, respectively (cf.\ \cite{to1}, Corollary $13$). It
suffices to argue that for such a Riemannian manifold the pair of
orthogonal distributions $\Delta_1=\spa\{e_1,e_2\}$ and
$\Delta_2=\spa\{e_3,e_4\}$ is invariant by conformal changes of the
metric. For that, let $W$ denote the Weyl curvature tensor of
$\Q_c^2\times \Sf^2$. Then
$$
\<W(e_i,e_j)e_l,e_k\>=(1+c)/3
$$
if $(i,j,k,l)\in\{(1,2,2,1),\, (2,1,1,2),\, (3,4,4,3,),\, (4,3,3,4)\}$,
$$
\<W(e_i,e_j)e_l,e_k\>=-(1+c)/3
$$
if $(i,j,k,l)\in \{(1,2,1,2),\, (2,1,2,1),\, (3,4,3,4,),\, (4,3,4,3)\}$, and
vanishes otherwise. 
In particular, if $X=\sum_{i=1}^4a_ie_i$ and
$Y=\sum_{j=1}^4b_je_j$ are orthonormal vectors, then
$$
\<W(X,Y)Y,X\> =\frac{1+c}{3}\left((a_1b_2-a_2b_1)^2
+(a_3b_4-a_4b_3)^2\right).
$$
Hence $\<W(X,Y)Y,X\>=0$ if and only if there exist $\lambda\in \R$, $Z\in
\spa\{e_1,e_2\}$ and \mbox{$U\in \spa\{e_3,e_4\}$} such that $X=\lambda
Z+U$ and $Y=Z-\lambda U$. In other words,
$$
\<W(X,Y)Y,X\>=0
$$
if and only if there exist vectors $S\in \spa\{e_1,e_2\}$ and $T\in
\spa\{e_3,e_4\}$ satisfying that $\spa\{X,Y\}=\spa\{S,T\}$. By the
conformal  invariance of $W$, the pair of orthogonal distributions
$\Delta_1$ and $\Delta_2$ is uniquely determined up to conformal changes
of the metric by the fact that $W(\sigma)=0$ for a two-plane $\sigma$ if
and only if $\sigma$ intersects both $\Delta_1$ and~$\Delta_2$.\qed

\begin{remark}\po\label{re:neq2} \textnormal{For $n=2$, the
proof of Lemma \ref{le:1}- $(ii)$ shows that $f$ and $\tilde{f}$ must have
a common mean curvature function $H=\alpha$. In particular, this implies
the well-known fact that if two surfaces $f\colon\,M^2\to \R^3$ and
$\tilde{f}\colon\,M^2\to \R^3$ induce conformal metrics on $M^2$ and
envelop a common sphere congruence, then the the latter is necessarily
their common central sphere congruence. Moreover, the proof can be easily
adapted to show that the same conclusion is true for a pair of immersions
$\tilde{f}, f\colon\,M^2\to \Q_c^3$ into any space form. }
\end{remark}

\begin{lemma}\po\label{le:2}
Under the assumptions of Lemma \ref{le:1} it holds that \be\label{B}
\nabla P(X,Y)+X\wedge Y ((P-I)\nabla\va)=0 \ee where $P=B\circ
T=e^{-\varphi}\tilde{B}$.
\end{lemma}

\proof Since $ \tilde{f}, f\colon\,M^n\to\R^{n+1}$, $n\geq 3$, induce
conformal metrics $\<\,,\,\>^\sim=e^{2\va}\<\,,\,\>$ on $M^n$, the
corresponding Levi-Civita connections $\tilde\nabla$ and $\nabla$ are
related by (\ref{eq:lcivita}). Using this we obtain that
$$
\begin{array}{l}
(\tilde{\nabla}_X\tilde{A})Y-(\tilde{\nabla}_Y\tilde{A})X \vspace{1ex} \\
=({\nabla}_X\tilde{A})Y -({\nabla}_Y\tilde{A})X
+(\tilde{A}Y) (\varphi)X-Y(\varphi)\tilde{A}X -(\tilde{A}X)(\varphi)Y
+X(\varphi)\tilde{A}Y.
\end{array}
$$
The Codazzi equation for $\tilde{f}$ then gives
$$
({\nabla}_X\tilde{A})Y-({\nabla}_Y\tilde{A})X =(\tilde{A}X)(\varphi)Y
-X(\varphi)\tilde{A}Y-(\tilde{A}Y)(\varphi)X +Y(\varphi)\tilde{A}X,
$$
which easily implies that
$$
\begin{array}{l}
({\nabla}_X\tilde{B})Y-({\nabla}_Y\tilde{B})X
+X(\alpha)Y-Y(\alpha)X\vspace{1ex}
\\
\hspace{20ex} =(\tilde{B}X)(\varphi)Y-X(\varphi)\tilde{B}Y
-(\tilde{B}Y)(\varphi)X+Y(\varphi)\tilde{B}X.
\end{array}
$$
Replacing $\tilde{B}=e^{-\varphi}(B\circ T)$ into the preceding equation
yields
$$
\begin{array}{l}
\<\nabla\varphi,B\circ T(X)\>Y-\<\nabla\varphi,B\circ T(Y)\>X \vspace{1ex}\\
\!=B((\nabla_XT)Y-(\nabla_YT)X)+(\nabla_XB)TY-(\nabla_YB)TX
+X(\alpha)Y-Y(\alpha)X,
\end{array}
$$
that is,
$$
(\nabla_XP)Y-(\nabla_YP)X =\<PX,\nabla\va\>Y-\<PY,\nabla\va\>X + X\wedge Y
(\nabla\va),
$$
which is equivalent to (\ref{B}). \qed

\begin{lemma}\po\label{le:3}
Under the assumptions of Lemma \ref{le:1}, suppose  
$\gamma:=\lambda-\alpha\neq 0$. Then, we have that

\begin{itemize}
\item[(i)] if the conclusion of Lemma \ref{le:1} holds
with $\epsilon=1$ then $\Delta^\perp$ is an umbilical distribution;

\item[(ii)] if the conclusion of Lemma \ref{le:1} holds
with $\epsilon=-1$ then $coker\, C$ has dimension $1$.
\end{itemize}
\end{lemma}

\proof \noindent \textit{Case $n\ge 4$.} Let $e_1,\ldots, e_n$ be an
orthonormal frame field such that $Be_1=\beta e_1$, $Be_2=-\beta e_2$ and
$Be_i=\gamma e_i$ for $i\geq 3$, where $\beta=\mu_1-\alpha=\alpha-\mu_2$.
Applying (\ref{B}) for $X=e_i$ and $Y=e_j$, $i\neq j\geq 3$, we obtain
that $\nabla\varphi\in \Delta^\perp$. Since $\det T|_{\Delta^\perp}=1$, we
may set
$$
Te_1=\cos\theta e_1+\sin\theta e_2,\,\,\,\,Te_2=-\sin\theta e_1+\cos\theta
e_2
$$
for some smooth function $\theta$. Since $A$ and $\tilde{A}$ can not be
simultaneously diagonalized at any point of $M^n$, we have that
$\cos\theta\neq 1$ everywhere.

The Codazzi equation for $f$ yields \be\label{cod0}
(\nabla_XB)Y+X(\alpha)Y=(\nabla_YB)X+Y(\alpha)X \ee for all tangent vector
fields $X,Y$. Applying (\ref{cod0}) for $X=e_1$ and $Y=e_2$ and taking the
$e_i$-component for $i\geq 3$ yields \be\label{cod1}
(\gamma+\beta)\omega_{i2}(e_1) =(\gamma-\beta)\omega_{i1}(e_2),\,\,\,i\geq
3. \ee Similarly, using (\ref{cod0}) for $X=e_1$ and $Y=e_i$, and then for
$X=e_2$ and $Y=e_i$, $i\geq 3$, and comparing the $e_1$-component of the
former equation with the $e_2$-component of the latter, we get
\be\label{cod2} (\gamma-\beta)\omega_{i1}(e_1)
+(\gamma+\beta)\omega_{i2}(e_2) =2e_i(\alpha),\,\,\,i\geq 3. \ee

Applying (\ref{B}) for $X=e_1$ and $Y=e_i$, $i\geq 3$, using (\ref{cod0})
and taking the $e_1$ and $e_2$-components yield, respectively,
\be\label{a}
\gamma(\epsilon-\cos\theta)\omega_{i1}(e_1)
+\beta\sin\theta\omega_{i2}(e_1)+\beta\sin\theta
e_i(\theta)+2\beta\sin\theta\omega_{21}(e_i)=(1-\cos\theta)e_i(\alpha)
\ee
and
\be\label{b}
\gamma(\epsilon-\cos\theta)\omega_{i2}(e_1)
+\beta\sin\theta\omega_{i1}(e_1)+ \beta\cos\theta
e_i(\theta)
-(\gamma+\beta)\sin\theta\omega_{i2}(e_2)=-\sin\theta e_i(\alpha).
\ee

Similarly, applying (\ref{B}) for $X=e_2$ and $Y=e_i$, $i\geq 3$, using
(\ref{cod0}) and taking the $e_1$ and $e_2$-components yield,
respectively,
\be\label{c}
\gamma(\epsilon-\cos\theta)\omega_{i1}(e_2)
+\beta\sin\theta\omega_{i2}(e_2)+\beta\cos\theta e_i(\theta)
(\gamma-\beta)\sin\theta\omega_{i1}(e_1)=\sin\theta e_i(\alpha)
\ee
and
\be\label{d}
\gamma(\epsilon-\cos\theta)\omega_{i2}(e_2)
+\beta\sin\theta\omega_{i1}(e_2)-\beta\sin\theta
e_i(\theta)+2\beta\sin\theta\omega_{12}(e_i)
=(1-\cos\theta)e_i(\alpha).
\ee
Subtracting (\ref{c}) from (\ref{b}) and adding 
(\ref{a}) and (\ref{d}) yields, respectively, \be\label{sub}
\gamma(\epsilon-\cos\theta)(\omega_{i2}(e_1)
-\omega_{i1}(e_2))+\beta\sin\theta(\omega_{i1}(e_1) -\omega_{i2}(e_2))=0
\ee 
and 
\be\label{add} 
\begin{array}{l}
\gamma(\epsilon-\cos\theta)(\omega_{i1}(e_1)
+\omega_{i2}(e_2))+\beta\sin\theta(\omega_{i2}(e_1) +\omega_{i1}(e_2))
=2(1-\cos\theta)e_i(\alpha). 
\end{array}
\ee

We now prove $(ii)$. By Lemma \ref{le:kerC}, if $\mbox{coker}\,C$ does not
have dimension~$1$ then there exists $S\in \mbox{coker}\,C$ such that
$C_S=a I$ for some nonzero real number~$a$. It follows from (\ref{cod2})
and (\ref{add}), respectively, that
$$
-a\gamma=S(\alpha)
$$
and
$$
a\gamma(\epsilon-\cos\theta)=(1-\cos\theta)S(\alpha),
$$
so we get a contradiction since $\gamma\neq 0$.

In order to prove $(i)$, we regard (\ref{cod1}), (\ref{cod2}),
(\ref{sub}) and (\ref{add})  as a system of linear equations in the
unknowns $\omega_{i1}(e_1)$, $\omega_{i2}(e_1)$, $\omega_{i1}(e_2)$ and
$\omega_{i2}(e_2)$. We obtain that its unique solution is
$\omega_{i1}(e_1)=\omega_{i2}(e_2)=e_i(\alpha)/\gamma$ and
$\omega_{i2}(e_1)=\omega_{i1}(e_2)=0$, which implies that $\Delta^\perp$
is
an umbilical distribution with mean curvature vector field $\eta=(1/\gamma)(\nabla\,\alpha)|_{\Delta}$.\vspace{1ex}\\
\textit{Case $n=3$.} Here we reorder the principal frame $e_1, e_2, e_3$
so that $e_1$ spans~$\Delta$ and
$$
Te_1=e_1,\,\,\,Te_2=\cos \theta e_2+\sin \theta
e_3,\,\,\,Te_3=-\sin \theta e_2 +\cos \theta e_3.
$$
We have
$$
Be_1=\gamma e_1,\,\,Be_2=-\beta e_2\,\,\,\mbox{and}\,\,\,Be_3=\beta e_3,
$$
with $\alpha-\mu_2=\beta=\mu_3-\alpha$.
 From (\ref{cod0}) we obtain
\be\label{coda0} e_2(\gamma+\alpha)=(\gamma
+\beta)\omega_{12}(e_1),\,\,\,\, e_3(\gamma+\alpha)=(\gamma
-\beta)\omega_{13}(e_1), \ee \be\label{coda1} \omega_{12}(e_2)
=\frac{e_1(\alpha-\beta)}{\gamma+\beta},\,\,\,\,\,\,\, \omega_{13}(e_3)
=\frac{e_1(\beta+\alpha)}{\gamma-\beta}, \ee and \be\label{coda2}
\omega_{13}(e_2)=\frac{-2\beta}{\gamma -\beta}\omega_{23}(e_1),
\,\,\,\,\,\,\,\, \omega_{12}(e_3)=\frac{-2\beta}{\gamma
+\beta}\omega_{23}(e_1). \ee Taking into account the first equation in
(\ref{coda0}), the $e_1$-component of (\ref{B}) for $X=e_1$ and $Y=e_2$
gives \be\label{aa} \cos\theta e_2(\va)-\sin\theta
e_3(\va)=(\cos\theta-1)\omega_{12}(e_1) -\sin\theta\omega_{13}(e_1).\ee
Using the first equations in (\ref{coda1}) and (\ref{coda2}), the $e_2$
and $e_3$-components of (\ref{B}) for $X=e_1$ and $Y=e_2$ yield,
respectively, \be\label{ab} \gamma e_1(\va) =\beta e_1(\theta)\sin\theta-
\frac{2\gamma\beta}{\gamma-\beta}\sin\theta\omega_{23}(e_1)
-\frac{(\cos\theta-1)}{\gamma+\beta}(\gamma e_1(\beta)+ \beta e_1(\alpha))
\ee and \be\label{ac} \beta\cos\theta
e_1(\theta)-\frac{2\gamma\beta}{\gamma-\beta}(\cos\theta-1)\omega_{23}(e_1)
+\frac{\sin\theta}{\gamma+\beta}(\gamma e_1(\beta)+\beta
e_1(\alpha))=0.\ee Similarly, taking the $e_1$-component of (\ref{B}) for
$X=e_1$ and $Y=e_3$ and using the second equation in (\ref{coda0}) give
\be\label{ba} \sin\theta e_2(\va)+\cos\theta
e_3(\va)=-\sin\theta\omega_{12}(e_1) -(\cos\theta-1)\omega_{13}(e_1). \ee
Using the second equations in (\ref{coda1}) and (\ref{coda2}), the
$e_3$-component of (\ref{B}) for $X=e_1$ and $Y=e_3$ yields \be\label{bb}
\gamma e_1(\va)=-\beta e_1(\theta)\sin\theta +\frac{2\gamma\beta}
{\gamma+\beta}\sin\theta\omega_{23}(e_1)
+\frac{(\cos\theta-1)}{\gamma-\beta}(\gamma e_1(\beta) +\beta
e_1(\alpha)). \ee Now, using all the equations in (\ref{coda1}) and
(\ref{coda2}) we obtain by taking the \mbox{$e_1$-component} of (\ref{B}) for
$X=e_2$ and $Y=e_3$ that \be\label{ca}
2\gamma\beta(\cos\theta-1)\omega_{23}(e_1) +\sin\theta(\gamma
e_1(\beta)+\beta e_1(\alpha))=0. \ee It follows from (\ref{ac}) and
(\ref{ca}) that \be\label{f1} \beta\cos\theta e_1(\theta) +
\frac{2\gamma\sin \theta}{\gamma^2-\beta^2}(\gamma e_1(\beta)+\beta
e_1(\alpha))=0. \ee On the other hand, using (\ref{ca}) we get from
(\ref{ab}) and (\ref{bb}) that \be\label{f2} \beta\sin\theta e_1(\theta) -
\frac{2\gamma\cos \theta}{\gamma^2-\beta^2}(\gamma e_1(\beta)+\beta
e_1(\alpha))=0.\ee Since $\beta\gamma\neq 0$, for $\gamma\neq 0$ by
assumption and $\beta\neq 0$ by Lemma \ref{le:1}-$(iii)$, we obtain from
(\ref{f1}) and (\ref{f2}) that \be\label{etheta} e_1(\theta)=0 \ee and \be
\label{umb1} \gamma e_1(\beta)+\beta e_1(\alpha)=0. \ee Then, it follows
from (\ref{coda1}) that \be\label{umb2} \omega_{12}(e_2)=\omega_{13}(e_3).
\ee In view of (\ref{umb1}), equations (\ref{ab}) and (\ref{bb}) reduce,
respectively, to
$$
\gamma e_1(\va)+\frac{2\gamma\beta}{\gamma
-\beta}\sin\theta\omega_{23}(e_1)=0
$$
and
$$
\gamma e_1(\va)
-\frac{2\gamma\beta}{\gamma+\beta}\sin\theta\omega_{23}(e_1)=0,
$$
which imply that \be\label{vaw23} e_1(\va)=0=\omega_{23}(e_1). \ee Then we
obtain from (\ref{coda2}) that
$$
\omega_{13}(e_2)=0=\omega_{12}(e_3).
$$
Together with (\ref{umb2}), this implies that 
$\Delta^\perp=\spa\{e_2,e_3\}$ is an umbilical 
distribution.\qed

\begin{proposition}\po\label{le:4}
Under the assumptions of Lemma \ref{le:1}, suppose further that its
conclusion holds with $\epsilon=1$ and that $\gamma:=\lambda-\alpha\neq
0$. Then both $f$ and $\tilde{f}$ are surface-like hypersurfaces.
\end{proposition}

\proof By Lemma \ref{le:3}, both $f$ and $\tilde{f}$ carry, respectively, principal
curvatures $\lambda$ and $\tilde{\lambda}$ of constant
multiplicity $n-2$ having a common eigenbundle $\Delta$ with the property
that $\Delta^\perp$ is an umbilical distribution.\vspace{1ex}

\noindent\textit{Case $n\ge 4$.} In this case the conclusion follows from
Theorem \ref{umb}. \vspace{1ex}

\noindent\textit{Case $n=3$.} We use the notations in the proof of Lemma
\ref{le:3}. Using that
$$
\tilde{A}=e^{-\va}BT+\alpha I
$$
and the Gauss equation for $\tilde{f}$, we obtain by taking the
$e_2$-component of (\ref{curvs}) for $X=e_1$ and $Y=Z=e_3$ that
\be\label{q1}Q(e_1,e_2)=0.\ee Similarly, the $e_3$-component of
(\ref{curvs}) for $X=e_1$ and $Y=Z=e_2$ gives \be\label{q2}Q(e_1,e_3)=0.
\ee In view of (\ref{vaw23}), equations (\ref{q1}) and (\ref{q2}) reduce
to
$$
e_1e_2(\va)=0=e_1e_3(\va),
$$
and (\ref{aa}), (\ref{ba}) and (\ref{etheta}) imply that \be\label{omegas}
e_1(\omega_{12}(e_1))=0=e_1(\omega_{13}(e_1)). \ee Then, using
(\ref{omegas}) and the second equality in (\ref{vaw23}), it follows that
the derivative of
$$
\nabla_{e_1}e_1=\omega_{12}(e_1)e_2+\omega_{13}(e_1)e_3
$$
with respect to $e_1$ has no $e_2$ and $e_3$-components. This means that
the integral curves of $e_1$ are circles in $M^3$.

On the other hand, it follows from  \cite{to1}, Lemma $12$  that also
$\Delta^\perp$ is a spherical distribution, that is, its mean curvature
vector field $\delta$ satisfies
$$
\<\nabla_Z\delta, e_1\>=0\;\;\mbox{for all}\;\; Z\in \Delta^\perp.
$$

We obtain from   \cite{to2}, Theorem 4.3 that $M^3$ is locally conformal
to a Riemannian product $I\times M^2$, the leaves of the product foliation
tangent to $I$ and $M^2$ corresponding, respectively, to the leaves of
$\Delta$ and $\Delta^\perp$. Then, Theorem 5 in \cite{to1} implies that
$f(M^3)$ is locally the image by a conformal transformation of Euclidean
space of one of the following: a cylinder $N^2\times\R$ over a surface
$N^2\subset\R^3$, a cone $CN^2$ over a surface $N^2\subset \Sf^3$, a
rotation hypersurface over a surface $N^2\subset \R_+^3$, a cylinder
$\gamma\times\R^2$ over a plane curve, a product $C\gamma\times\R$, where
$C\gamma$ is the cone over a spherical curve, or a rotation hypersurface
over a plane curve $\gamma\subset\R_+^2$. In the three last cases, the
hypersurface would have a principal curvature of multiplicity two with
$\Delta^\perp$ as eigenbundle, which is not possible by the assumption
that $A$ and $\tilde{A}$ can not be simultaneously diagonalized. Therefore
$f(M)$ and $\tilde{f}(M)$ must be (globally) open subsets of
hypersurfaces as in one of the first three cases.\qed

\section[The main lemma ]{The main lemma}

We now use the results of the previous section to show that a sphere
congruence in $\R^{n+1}$, $n\geq 3$, with conformal envelopes is
necessarily Ribaucour.

\begin{lemma}\po\label{prop:blaschke1}
Let $f,\tilde{f}\colon\,M^n\to\R^{n+1}$, $n\geq 3$, be a solution of
Blaschke's problem. Then the shape operators  $A$ and $\tilde{A}$ of $f$
and $\tilde{f}$, respectively, can be simultaneously diagonalized at any
point of $M^n$.
\end{lemma}

\proof As pointed out after Corollary \ref{cor:*}, condition $(ii)$ in
that result holds for $f$ and $\tilde{f}$, hence so does the conclusion of
Lemma \ref{le:1}. Since the set  of solutions of Blaschke's problem is
invariant under M\"obius transformations of $\R^{n+1}$,  by composing $f$
and $\tilde{f}$ with such a transformation we may assume, in view of
(\ref{pcurv}), that the function $\gamma$ defined in Lemma \ref{le:3} does
not vanish. For the same reason and bearing in mind (\ref{coker}), we may
suppose that the dimension of ${\mbox coker(C)}$ is not equal to $1$,
unless $\Delta^\perp$ is an umbilical distribution. It follows from Lemma
\ref{le:3} that $\Delta^\perp$ is indeed umbilical. By Theorem \ref{umb},
$f(M)$ and $\tilde{f}(M)$ are, up to (possibly different) M\"obius
transformations $\mathcal{I}_1$ and $\mathcal{I}_2$ of $\R^{n+1}$,
respectively, open subsets of one of the following:

\begin{itemize}
\item[(i)] cylinders $M^2\times\R^{n-2}$ and $\tilde{M}^2\times
\R^{n-2}$ over surfaces $M^2, \tilde{M}^2\subset\R^3$, respectively;

\item[(ii)] cylinders $CM^2\times\R^{n-3}$ and $C\tilde{M}^2\times
\R^{n-3}$, respectively, where $CM^2\subset\R^4$ denotes the cone over
$M^2\subset \mathbb{S}^3$;

\item[(iii)] rotation hypersurfaces over surfaces $M^2$ and $\tilde{M}^2$ contained
in $\R^3_+$, respectively.
\end{itemize}

 We now make the following key observation.

\begin{lemma}\po\label{le:obs} If $f$ is as in $(ii)$,
then $M^n$ is conformal to $M^2\times \Hy^{n-2}$ with $M^2$ endowed with
the metric induced from $\Sf^3$. If $f$ is as in $(iii)$,  then $M^n$ is
conformal to $M^2\times\Sf^{n-2}$ with $M^2$ endowed with the metric induced
from the metric of constant sectional curvature $-1$ on  $\R_+^3$ regarded
as the half-space model of $\Hy^3$.
\end{lemma}

\proof  An $f$ as in $(ii)$ can be parameterized by
$$
f=(t_1, \ldots, t_{n-3}, t_{n-2}g_1, \ldots, t_{n-2}g_4),
$$
where $g=(g_1, \ldots, g_4)$ parameterizes $M^2\subset \Sf^3$. Then the
metric induced by $f$ is
$$
\<\,,\,\>=\<\,,\,\>_{\R^{n-2}}+t_{n-2}^2\<\,,\,\>_{M^2}
=t_{n-2}^2(\<\,,\,\>_{\Hy^{n-2}}+\<\,,\,\>_{M^2}),
$$
where $\<\,,\,\>_{\Hy^{n-2}}$ is the metric of constant sectional
curvature $-1$ on $\R_+^{n-2}$ regarded as the half-space model of
$\Hy^{n-2}$.

 Any $f$  as in $(iii)$ can be parameterized by
$f=(g_1, g_2, g_3\phi)$, where $g=(g_1,g_2,g_3)$ parameterizes
$M^2\subset\R^3_+$ and $\phi$ parameterizes $\Sf^{n-2}\subset\R^{n-1}$.
Thus the metric induced by $f$ is
$$
\<\,,\,\>=\<\,,\,\>_{M^2}+g_3^2\<\,,\,\>_{\Sf^{n-2}}
=g_3^2(\<\,,\,\>^\sim_{M^2}+\<\,,\,\>_{\Sf^{n-2}}),
$$
where $\<\,,\,\>_{M^2}$ denotes the metric  on $M^2$ induced by the
Euclidean metric on $\R^3_+$ and  $\<\,,\,\>^\sim_{M^2}$ the metric on
$M^2$ induced by the hyperbolic metric of constant sectional curvature
$-1$ on  $\R^3_+$ regarded as the half-space model 
of~$\Hy^{3}$.\qed

\begin{remark}\po \label{re:metric}
\textnormal{If $f\colon\,M^n\to\R^{n+1}$ is either a cylinder
$C(\gamma)\times\R^{n-2}$, where $C(\gamma)$ is the cone over a curve
$\gamma\colon\,I\to\Sf^2\subset\R^3$, or a rotation hypersurface over a
curve $\gamma\colon\,I\to\R_+^2$, then the proof of Lemma \ref{le:obs}
shows that the metric induced by $f$ is conformal to $ds^2+d\sigma^2$,
where $s$ denotes the arc-length function of $\gamma$, regarded as a curve
in the half-space model of the hyperbolic plane in the last case, and
$d\sigma^2$ is the metric of either hyperbolic space $\Hy^{n-1}$ or
the sphere $\Sf^{n-1}$ respectively. }
\end{remark}

\begin{lemma}\po\label{le:triv} Let
$\<\,,\,\>=\pi_1^*\<\,,\,\>_1+\pi_2^*\<\,,\,\>_2$ and
$\<\,,\,\>^\sim=\pi_1^*\<\,,\,\>_1^\sim +\pi_2^*\<\,,\,\>^\sim_2$ be
product metrics on a product manifold $M=M_1\times M_2$, where $\pi_i$
denotes the projection of $M_1\times M_2$ onto $M_i$ for $i=1,2$.  If
$\<\,,\,\>^\sim=\psi^2\<\,,\,\>$ for some $\psi\in\mathcal{C}^\infty(M)$,
then  $\psi$ is a constant $k\in\R$ and $\<\,,\,\>_i^\sim=k^2\<\,,\,\>_i$
for $i=1,2$.
\end{lemma}

\proof Given $i\in\{1,2\}$, let $X_i\in TM_i$ be a local unit vector field
with respect to $\<\,,\,\>_i$. Let $\tilde{X}_i$ be the lift of $X_i$ to
$M$. Then $\<X_i,X_i\>^\sim_i\circ\pi_i
=\<\tilde{X}_i,\tilde{X}_i\>^\sim=\psi^2,$ and  the conclusion
follows.\qed\vspace{1ex}

Going back to the proof of Lemma \ref{prop:blaschke1}, we obtain from
Lemmas \ref{le:obs} and \ref{le:triv} that $f$ and $\tilde{f}$ differ by
M\"obius transformations  $\mathcal{I}_1$ and $\mathcal{I}_2$,
respectively, from surface-like hypersurfaces that are of the same type,
with the corresponding surfaces $M^2$ and $\tilde{M}^2$ being isometric as
surfaces as in either $\R^3$, $\Sf^3$ or $\Hy^3$, respectively.

 Using again the invariance of the set of
solutions of Blaschke's problem under M\"obius transformations of
$\R^{n+1}$, we may assume that $\mathcal{I}_1$ is the identity map of
$\R^{n+1}$. Since the shape operators of $f$ and $\tilde{f}$ satisfy
$\mbox{trace }(A|_{\Delta^\perp})=\mbox{trace
}(\tilde{A}|_{\Delta^\perp})$, as follows from Lemma~\ref{le:1}, we obtain
that $\mbox{trace }(\tilde{A}|_{\Delta^\perp})$ is constant along
$\Delta$. Taking (\ref{pcurv}) into account once more, this easily implies
that $\mathcal{I}_2$ must be a similarity in cases $(i)$ and $(ii)$, and a
composition of a similarity with an inversion with respect to a
hypersphere centered at a point of the rotation axis in case $(iii)$. In
either case $\mathcal{I}_2\circ\tilde{f}$ is still a hypersurface as in
$(i)$, $(ii)$ or $(iii)$, respectively.

In summary, we can assume that $f$  and $\tilde{f}$ are both as in $(i)$,
$(ii)$ or  $(iii)$. We now use that \be\label{env2}
f+RN=\tilde{f}+R\tilde{N}, \ee and that
$$
2/R=\mbox{trace }(A|_{\Delta^\perp}) =\mbox{trace
}(\tilde{A}|_{\Delta^\perp}),
$$
by Lemma \ref{le:1}. Suppose first that $f$  and $\tilde{f}$ are (open
subsets of) cylinders $M^2\times\R^{n-2}$ and $\tilde{M}^2\times\R^{n-2}$,
respectively. Then $R$ is constant along $\Delta$ and we can write
(\ref{env2}) as

\be\label{cyl} g(q)+t+R(q)N(q)=\tilde{g}(q)+\tilde{t}+R(q)\tilde{N}(q),
\ee where $g$ and $\tilde{g}$ denote the position vectors of $M^2$ and
$\tilde{M}^2$, respectively, and $t,\tilde{t}$ parameterize the rulings of
$f$ and $\tilde{f}$, respectively. Differentiating (\ref{cyl}) with
respect to a vector $T\in \Delta$ implies that $f$ and $\tilde{f}$ are
cylinders with respect to the same orthogonal decomposition
$\R^{n+1}=\R^3\times\R^{n-2}$. The same argument shows that also in case
$(ii)$ the hypersurfaces  $f$ and $\tilde{f}$ are cylinders with respect
to the same orthogonal decomposition $\R^{n+1}=\R^4\times\R^{n-3}$. We now
show that  $CM^2$ and $C\tilde{M}^2$ are cones over surfaces $M^2$,
$\tilde{M}^2$ in the same hypersphere of $\R^4$. In fact, if $g$ and
$\tilde{g}$ denote the position vectors of $M^2$ and $\tilde{M}^2$,
respectively, then, disregarding the common components in $\R^{n-3}$, we
can now write (\ref{env2}) as
$$
P_0+t(g(q)-P_0)+t(1/H(q))N(q) =\tilde{P}_0+t(\tilde{g}(q)
-\tilde{P}_0)+t(1/H(q))\tilde{N}(q),
$$
where $P_0$ and $\tilde{P}_0$ are the vertices of $CM^2$ and
$C\tilde{M}^2$, respectively, $H$ is the common mean curvature function of
$M^2$ and $\tilde{M}^2$, and $N$ and $\tilde{N}$ are unit normal vector
fields to $M^2$ and $\tilde{M}^2$, respectively. Letting $t$ go to $0$
yields $P_0=\tilde{P}_0$, as asserted.

 In case $(iii)$,  we claim that the rotation hypersurfaces $f$ and $\tilde{f}$ must
 have the same ``axis". In fact, in this case (\ref{env2}) can be written as
\be\label{eq:rot}\begin{array}{l} P_0+(g_1(q),g_2(q),g_3(q)\phi(t))+
R(q)(N_1(q),N_2(q),N_3(q)\phi(t))\vspace{1ex}
\\\hspace*{4ex}=
\tilde{P}_0+(\tilde{g}_1(q),\tilde{g}_2(q), \tilde{g}_3(q)\phi(t))+
R(q)(\tilde{N}_1(q),\tilde{N}_2(q), \tilde{N}_3(q)\phi(t)),\end{array} \ee
where $(g_1,g_2,g_3)$ and $(\tilde{g}_1,\tilde{g}_2,\tilde{g}_3)$
parameterize $M^2$ and $\tilde{M}^2$, $N$ and $\tilde{N}$ are unit normal
vector fields to $M^2$ and $\tilde{M}^2$, respectively, and $\phi$
parameterizes the unit sphere in $\R^{n-2}$. Coordinates in different
sides of (\ref{eq:rot}) are possibly with respect to different orthonormal
bases of $\R^{n+1}$. Differentiating  (\ref{eq:rot}) with respect to $t_i$
yields
$$
(f_3(q)+R(q)N_3(q))\frac{\d \phi}{\d t_i}
=(\tilde{f}_3(q)+R(q)\tilde{N}_3(q))\frac{\d\phi}{\d t_i},
$$
which implies that the rotation axes coincide up to translation. In
particular, we may choose a common orthonormal basis of $\R^{n+1}$ to
parameterize $f$ and $\tilde{f}$ as in (\ref{eq:rot}). Then we obtain
$$
(P_0)_i-(\tilde{P}_0)_i=\phi_{i+2}(t)(-g_3(q)-R(q)N_3(q)
+\tilde{g}_3(q)+R(q)\tilde{N}_3(q)),
$$
which implies that $(P_0)_i=(\tilde{P}_0)_i$ for $i=3,\ldots, n+1$, and
proves our claim.

Now,  intersecting the spheres of the congruence enveloped by $f$ and
$\tilde{f}$ with either the three-dimensional subspace $\R^3$ orthogonal
to their common Euclidean factor $\R^{n-2}$ in case $(i)$, the hypersphere
$\mathbb{S}^3\subset\R^4$  containing $M^2$ and $\tilde{M}^2$ in case
$(ii)$, or the affine subspace $\R^3$ of $\R^{n+1}$ containing the
profiles $M^2$ and $\tilde{M}^2$ in case $(iii)$, gives a sphere
congruence in $\R^3$, $\Sf^3$ and $\Hy^3$, respectively,  which is
enveloped by $M^2$ and  $\tilde{M}^2$.  Moreover, since the shape
operators of $f$ and $\tilde{f}$ can not be simultaneously diagonalized at
any point, the same holds for the shape operators of $M^2$ and
$\tilde{M}^2$.\vspace{1ex}

The proof  of Lemma \ref{prop:blaschke1} is now completed by the following
result.

\begin{proposition}\po\label{prop:bonnet1}
If $f_1, {f}_2\colon\,M^2\to\mathbb{Q}_c^3$ are isometric immersions whose
shape operators can not be simultaneously diagonalized at any point of
$M^2$ then they can not envelop a common sphere congruence.
\end{proposition}

\proof Assuming the contrary, the common sphere congruence enveloped by
$f$ and $\tilde{f}$ is necessarily their common central sphere congruence
(see Remark \ref{re:neq2}). In particular, $f_1$ and $f_2$ have the same
mean curvature function, and hence $(f_1,f_2)$ is a solution of Bonnet's
problem in $\Q_c^3$. We argue separately for each of the three types of
solutions of that problem, as discussed (for $\R^3$) at the beginning of
the introduction.

Assume first that $f_1, {f}_2\colon\,M^2\to \mathbb{Q}_c^3$ are isometric
immersions with the same {\em constant\/} mean curvature function. Using
that the radius function  of the sphere congruence enveloped by $f$ and
$\tilde{f}$  is constant we obtain from (\ref{der}) that
$$
f_*BX=\tilde{f}_*\tilde{B}X
$$
for every $X\in TM^2$, where $B=A-HI$ and $\tilde{B}=\tilde{A}-HI$. Since
$\tilde{B}$ is invertible by  Lemma \ref{le:1}-$(iii)$, the preceding
equation can be rewritten as
$$
\tilde{f}_*X=f_*\Phi X
$$
for every $X\in TM^2$, where $\Phi=B\circ \tilde{B}^{-1}$. Regarding
$\omega=\tilde{f}_*$ as a one-form in $M^2$ with values in either $\R^3$,
$\R^4$ or $\Les^4$, according as $c=0,1$ or $-1$, respectively, we obtain
by taking the normal component to $M^2$ in the equation $d\omega=0$ that
\be\label{eq:comut} \Phi^tA_\xi=A_\xi \Phi \ee for every normal vector
field $\xi$ to $M^2$, regarded as a surface in $\R^3$, $\R^4$ or $\Les^4$,
respectively (see  \cite{dt2}, Proposition $1$). Moreover, the fact that
$f$ and $\tilde{f}$ are isometric implies that  $\Phi$ is an orthogonal
tensor. Then we obtain from (\ref{eq:comut}) that $\mbox{trace} A_\xi=0$
for every normal vector $\xi$, a contradiction.

If $f$ and $\tilde{f}$ are  Bonnet surfaces admitting a one-parameter
family of isometric deformations preserving the Gauss map then they are
isothermic surfaces (cf.\ \cite{bob}). But isothermic surfaces  are
characterized as the only surfaces whose central sphere congruence is
Ribaucour (see \cite{hj}, Lemma $3.6.1$). Thus, this case is also ruled
out.

Finally,  suppose that  $f_1,f_2\colon\,M^2\to\mathbb{Q}_c^3$ form   a
Bonnet pair of surfaces with mean curvature function $H$ with respect to
unit normal vector fields $\eta_j$, $1\le j\le 2$. For simplicity we take
$c=0$ or $c=\pm 1$. We have from \cite{bob} and \cite{te} that in
conformal coordinates
$$
ds^2= e^u(dx^2+ dy^2)
$$
their second fundamental forms are given by \be  \label{sff} \left\{
\begin{array}{l}
A_{\eta_j}X = (H+he^{-u})X + \epsilon\, ke^{-u}Y \vspace*{1.5ex} \\
A_{\eta_j}Y= \epsilon\, ke^{-u}X + (H-he^{-u})Y,
\end{array}
\right. \ee where $X=\partial/\partial x$, $Y=\partial/\partial y$ are the
coordinate tangent vector fields, $\epsilon =1$ if $j=1$ and $\epsilon=-1$
if $j=2$. Moreover, $k\neq 0$ is a constant and $h\in
\mathcal{C}^\infty(M)$ satisfies the Codazzi equations \be  \label{1}
\left\{
\begin{array}{l}
H_xe^u-h_x=0 \vspace*{1.5ex} \\
H_ye^u+h_y=0.
\end{array}
\right. \ee The integrability condition of (\ref{1}) is \be  \label{int}
2H_{xy}+H_xu_y+H_y+H_x=0. \ee We first consider the case $c=0$. We show
that the surfaces
$$
F^j=f_j+\frac{1}{H}\eta_j,\;\;\;\; 1\le j\le 2,
$$
are isometric but never isometrically congruent. The coordinate vector
fields of $F^j$ are given by
$$
\left\{
\begin{array}{l}
F^j_*X= {\displaystyle \frac{-e^{-u}}{H}(\epsilon\, kY+hX)
+ \left(\frac{1}{H}\right)_x\eta_j}\vspace*{1.5ex} \\
F^j_*Y= {\displaystyle \frac{-e^{-u}}{H}(\epsilon\, kX-hY) +
\left(\frac{1}{H}\right)_y\eta_j.}
\end{array}
\right.
$$
Then,
$$
N_j=(\epsilon\,kH_y+hH_x)X + (\epsilon\,kH_x-hH_y)Y - H(k^2+h^2)\eta_j
$$
is a normal vector field to $F_j$. Thus, we have
$$
\left\{
\begin{array}{l}
\|F^1_*X\|={\displaystyle \frac{e^{-u}}{H^2}(h^2+k^2) + \frac{H^2_x}{H^4}=\|F^2_*X\|}\vspace*{1.5ex} \\
\|F^1_*Y\|={\displaystyle \frac{e^{-u}}{H^2}(h^2+k^2) +
\frac{H^2_y}{H^4}=\|F^2_*Y\|}
\vspace*{1.5ex} \\
\<F^1_*X,F^1_*Y\>={\displaystyle \frac{H_xH_y}{H^4} =\<F^2_*X,F^2_*Y\>,}
\end{array}
\right.
$$
and
$$
\|N_j\|^2=(h^2+k^2)(H^2+e^u(H_x^2+H_y^2)).
$$
Therefore, the surfaces $F^1,F^2$ are isometric and $\|N_1\|=\|N_2\|$.

Now, a long but straightforward computation using (\ref{sff}), (\ref{1})
and (\ref{int}) shows that the second fundamental forms $B^j_N$ of $F_j$,
$j=1,2$, satisfy that
$$
\<B^1_NX,F^1_*X\>=\<B^2_NX,F^2_*X\>,\;\;
\<B^1_NX,F^1_*Y\>=\<B^2_NX,F^2_*Y\>
$$
and
$$
\<B^1_NX,F^1_*Y\>=\<B^2_NX,F^2_*Y\>
$$
are independent of $\epsilon$. On the other hand,
$$
\<B^1_NX,F^1_*Y\>= \epsilon\,\frac{k}{H}\left(e^u(H_ x^2+H_y^2) +H^2(k^2+
h^2)\right)=\<B^2_NX,F^2_*Y\>
$$
never vanishes. Thus the surfaces cannot be isometrically
congruent.\vspace{1ex}

We now consider the case $c\neq 0$. We take $\mathbb{S}^3\subset\R^4$ and
$\mathbb{H}^3\subset\mathbb{L}^4$, and thus $\<f_j, f_j\>=c$. As before,
we show that the surfaces
$$
F^j=Cf_j+S\eta_j,\;\;\;\; 1\le j\le 2,
$$
are isometric but never isometrically congruent. Here
$$
\left\{
\begin{array}{l}
C=\cos R,\;\;\; S=\sin R\;\;\;\;\;\mbox{where}\;\;\, \cot R=H\;\;\;\;\;\;\mbox{if}\; c=1, \vspace*{1.5ex} \\
C=\cosh R,\;S=\sinh R\;\;\;\mbox{where}\;\; \coth R=H\;\;\;\;\mbox{if}\;\;
c=-1.\end{array} \right.
$$
Using (\ref{sff}) we easily obtain that the coordinate vector fields of
$F^j$ are
$$
\left\{
\begin{array}{l}
F^j_*X= R_x(C\eta_j-cSf_j) -She^{-u}X -\epsilon kSe^{-u}Y, \vspace*{1.5ex}
\\
F^j_*Y= R_y(C\eta_j-cSf_j) + She^{-u}Y -\epsilon kSe^{-u}X,
\end{array}
\right.
$$
where $R_x=-S^2H_x$ and $R_y=-S^2H_y$. Then,
$$
N_j=S(h^2+k^2)(C\eta_j-cSf_j)+ (hR_x+\epsilon kR_y)X - (hR_y-\epsilon
kR_x)Y
$$
is a normal vector field to $F_j$. Then, we have
$$
\left\{
\begin{array}{l}
\|F^1_*X\|=S^2e^{-u}(h^2+k^2) + R^2_x =\|F^2_*X\| \vspace*{1.5ex} \\
\|F^1_*Y\|=S^2e^{-u}(h^2+k^2) + R^2_y =\|F^2_*Y\| \vspace*{1.5ex} \\
\<F^1_*X,F^1_*Y\>= R_xR_y=\<F^2_*X,F^2_*Y\>,
\end{array}\right.
$$
and
$$
\|N_j\|^2=(h^2+k^2)(S^2(h^2+k^2)+e^u(R_x^2+R_y^2)).
$$
Therefore, the surfaces $F^1,F^2$ are isometric and $\|N_1\|=\|N_2\|$.

As before, the second fundamental forms $B^j_N$ of $F_j$, $j=1,2$, satisfy
that
$$
\<B^1_NX,F^1_*X\>=\<B^2_NX,F^2_*X\>,\;\;
\<B^1_NY,F^1_*Y\>=\<B^2_NY,F^2_*Y\>
$$
and
$$
\<B^1_NY,F^1_*Y\>=\<B^2_NY,F^2_*Y\>
$$
are independent of $\epsilon$. On the other hand,
$$
\<B^1_NX,F^1_*Y\>= -\epsilon\,\frac{k}{S}\left(e^u(H_ x^2+H_y^2)
+S^2C^2(k^2+ h^2)\right)=\<B^2_NX,F^2_*Y\>
$$
never vanishes, hence the surfaces cannot be isometrically congruent.\qed

\begin{remark}\po \textnormal{In Lemma \ref{prop:blaschke1} one does not need to
assume regularity of the enveloped sphere congruence.}
\end{remark}

\section[Proof of Theorem
\ref{thm:blaschke}]{Proof of Theorem \ref{thm:blaschke}}

 By Lemma \ref{prop:blaschke1}, the shape operators of $f$ and $\tilde{f}$ can be
simultaneously diagonalized at any point of $M^n$.  Moreover,  by
assumption the sphere congruence enveloped by $f$ and $\tilde{f}$ is
regular, thus the inverse of its radius function is nowhere a principal
curvature of either hypersurface  by Proposition~\ref{prop:indmetrics}.
Therefore $f$ and $\tilde{f}$ are Ribaucour transforms  one of each other,
as defined in \cite{dt2}. Furthermore, since they induce conformal metrics
on $M^n$, they are in fact \emph{Darboux} transforms one of each other in
the sense of \cite{to1}. Isometric immersions $f\colon\,M^n\to\R^{N}$,
$n\geq 3$, that admit Darboux transforms have been classified in
\cite{to1}, Theorem $20$. In order to state that result we first recall
some preliminary facts.

  Let $f\colon\,M^n\to\R^{N}$ be an isometric immersion
of a simply-connected Riemannian manifold with second fundamental form $\alpha\colon\,TM\times TM\to T_f^\perp M$. We have from
\cite{dt2}, Theorem $17$ that if $\tilde{f}\colon\,M^n\to\R^{N}$ is a
Ribaucour transform of $f$ then there exist $\va\in \mathcal{C}^\infty(M)$
and $\beta\in T_f^\perp M$ satisfying \be\label{eq:gnorm}
\alpha(X,\nabla\va)+\nabla_X^\perp\beta=0\;\; \mbox{for all}\;\; X\in TM
\ee such that \be\label{eq:rb} \tilde{f}= f - 2\nu\varphi \Fes, \ee where
$\Fes=df(\grad\,\va)+\beta$ and $\nu^{-1}=\<\Fes,\Fes\>$. Therefore
$\tilde{f}$ is completely determined by $(\va,\beta)$, or equivalently, by
$\va$ and $\Fes$. We denote \mbox{$\tilde{f}={\cal R}_{\va,\beta}(f)$.}
Moreover, we have that
$$
{\cal S}_{\va,\beta}:=\hess\va-A_\beta
$$
is a Codazzi tensor on $M^n$ such that
$$
\alpha({\cal S}_{\va,\beta}X,Y) =\alpha(X,{\cal S}_{\va,\beta}Y)\;\;
\mbox{for all}\;\; X,Y\in TM
$$
and
$$
d\Fes=df\circ {\cal S}_{\va,\beta}.
$$
Conversely, given $(\va,\beta)$ satisfying (\ref{eq:gnorm}) on an open
subset $U\subset M^n$ where
$$
D:=I- 2\nu\varphi{\cal S}_{\va,\beta}
$$
is invertible, then $\tilde{f}$ given by (\ref{eq:rb}) defines a \rtf of
$f|_U$, and the induced metrics of $f$ and $\tilde{f}$ are related by
\be\label{eq:met} \<X,Y\>^{\sim}=\<DX,DY\>. \ee

It follows from (\ref{eq:met}) and the symmetry of $D$ that if $f$ and
$\tilde{f}$ induce conformal metrics on $M^n$ then $D^2=r^2I$ for some
$r\in \mathcal{C}^\infty(M)$. Therefore, either  $D=\pm r I$ or $TM$
splits  orthogonally as $TM=E_+\oplus E_-$, where $E_+$ and $E_-$ are the
eigenbundles of $D$ correspondent to the eigenvalues $r$ and $-r$,
respectively. In the first case, it follows from the results in \cite{dt2}
that there exists an inversion $I$ in $\R^{N}$ such that
$L'(\tilde{f})=I(L(f))$, where $L$ and $L'$ are compositions of a
homothety and a translation.  The immersion  $\tilde{f}$ is said to be a
Darboux transform  of $f$ if the second possibility holds, in which case
$E_+$ and $E_-$ are also the eigenbundles of ${\cal S}_{\va, \beta}$ correspondent to
its distinct eigenvalues $\lambda=h(1-r)$ and $\mu=h(1+r)$, respectively,
where $h^{-1}=2\nu\va$. Thus, $\tilde{f}$ is a Darboux transform  of $f$
if and only if the associated Codazzi tensor ${\cal S}_{\va, \beta}$ has exactly two
distinct eigenvalues $\lambda, \mu$ everywhere satisfying
\be\label{eq:darboux} (\lambda+\mu)\va=\nu^{-1} =\<\Fes,\Fes\>. \ee The
classification of isometric immersions  $f\colon\,M^n\to\R^{N}$, $n\geq
3$, that admit Darboux transforms is as follows.

\begin{theorem}\po\label{prop:comb2}  Let $f\colon\,M^n\to\R^{N}$, $n\geq 3$, be an isometric
immersion that admits a Darboux transform $\tilde{f}={\cal
R}_{\va,\beta}(f)\colon\,M^n\to\R^{N}$. Then there exist locally a
product representation $\psi\colon\,M_1\times M_2\to M^n$ of $(E_+, E_-)$,
a homothety $H$ and an inversion $I$ in $\R^{N}$ such that one of the
following holds:
\begin{itemize}
\item[{\em (i)}] $\psi$ is a conformal diffeomorphism with respect to a Riemannian product metric
on $M_1\times M_2$ and
$$
f\circ \psi=H\circ  I\circ g,
$$
where $g=g_1\times g_2\colon\,M_1\times
M_2\to\R^{N_1}\times\R^{N_2}=\R^{N}$ is an extrinsic product of isometric
immersions. Moreover,  there exists $i\in \{1,2\}$ such that either $M_i$
is one-dimensional or $g_i(M_i)$ is contained in some sphere
$\Sf^{N_i-1}(P_i;r_i)\subset\R^{N_i}$.
\item[{\em (ii)}] $\psi$ is a conformal diffeomorphism
with respect to a warped product metric on \mbox{$M_1\times M_2$}  and
$$
f\circ \psi=H\circ I\circ\Phi\,\circ(g_1\times g_2),
$$
where $\Phi\colon\,\R_+^{m}\times_{\sigma}\Sf^{N-m}(1)\to\!\R^N$ denotes
an isometry and $g_1\colon M_1\to\R_+^{m}$ and \mbox{$g_2\colon
M_2\to\Sf^{N-m}(1)$} are isometric immersions.
\end{itemize}
Conversely, any such \ii admits a Darboux transform.
\end{theorem}

By a  product representation $\psi\colon\,M_1\times M_2\to M^n$ of the
orthogonal net $(E_+, E_-)$ we mean a diffeomorphism that maps the leaves
of the product foliation of $M_1\times M_2$ induced by $M_1$
(respectively, $M_2$) onto the leaves of $E_+$ (respectively, $E_-$).  In
part $(ii)$ the isometry $\Phi$ and the isometric immersions $g_1$ and
$g_2$ may be taken so that $g_2(M_2)$ is not contained in any hypersphere
of $\Sf^{N-m}(1)$.

The preceding definition of a Darboux transform
$\tilde{f}\colon\,M^n\to\R^{N}$ of an isometric immersion
$f\colon\,M^n\to\R^{N}$  does not rule out the possibility that
$\tilde{f}$ be conformally congruent to $f$. In fact, we now prove that in
either case of Theorem \ref{prop:comb2} this  always occurs whenever both
factors  $M_1$ and $M_2$ have dimension greater than one.

\begin{proposition}\po\label{prop:comb3} Let $f\colon\,M^n\to\R^{N}$ be as in $(i)$ or $(ii)$
 of Theorem~\ref{prop:comb2}. If both   $M_1$ and $M_2$
have dimension greater than one then any Darboux transform
$\tilde{f}\colon\,M^n\to\R^{N}$ of $f$ is conformally congruent to it.
\end{proposition}

For the proof of Proposition \ref{prop:comb3} we will need the following
fact on Codazzi tensors on warped products.

\begin{lemma}\po\label{le:codtensor2} Let $M^n=M_1\times_{\rho} M_2$ be a warped
product, let $(E_1,E_2)$ be its product net and let ${\cal
S}=\lambda\Pi_1+\mu\Pi_2$ be a Codazzi tensor on $M^n$ with $\lambda\neq
\mu$ everywhere, where $\Pi_i$ denotes orthogonal projection of $TM$ onto
$E_i$ for $i=1, 2$. If both factors have dimension greater than one then
$\lambda=A\in\R$ and $\mu=B(\rho\circ \pi_1)^{-1}+A$ for some $B\neq 0$.
\end{lemma}

\proof Since $E_1$ and $E_2$ are the eigenbundles of the Codazzi tensor
${\cal S}$, both $E_1$ and $E_2$ are umbilical  with mean curvature
normals given, respectively,
 (see \cite{re} or \cite{to2}, Proposition 5.1) by
\be\label{eq:etazeta} \eta=(\lambda-\mu)^{-1}(\nabla
\lambda)_{E_2}\;\;\;\mbox{and}\;\;\; \zeta=(\mu-\lambda)^{-1}(\nabla
\mu)_{E_1}. \ee Here, writing a vector subbundle as a subscript of a
vector field indicates taking the orthogonal projection of the vector
field onto that subbundle.

On the other hand, since $(E_1,E_2)$ is the product net of a warped
product with warping function $\rho$ we have (see \cite{mrs}, Proposition
$2$) \be\label{eq:etazeta2} \eta=0\;\;\;\mbox{and}\;\;\; \zeta=-\nabla
\log\circ \rho\circ \pi_1. \ee It follows from (\ref{eq:etazeta}) and
(\ref{eq:etazeta2}) that there exists $\tilde{\lambda}\in
\mathcal{C}^\infty(M_1)$ such that $\lambda=\tilde{\lambda}\circ \pi_1$.
Since any eigenvalue of a Codazzi tensor is constant along its eigenbundle
whenever the latter has rank greater than one (see \cite{to2}, Proposition
$5.1$), we obtain that $\tilde{\lambda}=A$ for some $A\in\R$ and that
$\mu$ is constant along $M_2$. Hence $\nabla \mu =(A-\mu)\nabla \log\circ
\rho\circ \pi_1$ by the second equation in (\ref{eq:etazeta}). This
implies  that $\nabla (\mu(\rho\circ\pi_1))=A\nabla \rho\circ \pi_1$, and
the conclusion follows.\qed

\begin{corollary}\po\label{le:codtensor1}
Let $M^n=M_1\times M_2$ be a Riemannian product, let $(E_1,E_2)$ be its
product net and let ${\cal S}=\lambda\Pi_1+\mu\Pi_2$ be a Codazzi tensor
on $M^n$ with $\lambda\neq \mu$ everywhere.  If both factors have
dimension greater than one then both  $\lambda$ and $\mu$ are
constants.\end{corollary}

\noindent {\em Proof of Proposition \ref{prop:comb3}:\/} We first consider
$f$ as in case $(i)$ of Theorem \ref{prop:comb2}. It suffices to prove
that under the assumption on the dimensions of $M_1$ and $M_2$ any Darboux
transform
$$
\tilde{g}=g-2\va\nu\Fes\colon\,M_1\times M_2\to\R^{N}
$$
of an extrinsic product $g=g_1\times g_2\colon\,M_1\times
M_2\to\R^{N_1}\times\R^{N_2}=\R^{N}$ of isometric immersions, such that
the eigenbundle net of the associated Codazzi tensor ${\cal
S}_{{\va},{\beta}}=\hess \va-A_\beta$ is the product net of $M_1\times
M_2$,  is conformally congruent to $g$.

Since both $M_1$ and $M_2$  have dimension greater than one, it follows
from Corollary~\ref{le:codtensor1} that  the eigenvalues of ${\cal
S}_{{\va},{\beta}}$ are constant, say, $a_1, a_2\in\R$. Integrating
$d\Fes=dg\circ{\cal S}_{{\va},{\beta}}$ gives
$$
\Fes=(a_1(g_1-P_1),a_2(g_2-P_2))
$$
for some $P_1\in\R^{N_1}$ and  $P_2\in\R^{N_2}$. Using that
${\va}(a_1+a_2)={{\nu}}^{-1}= \<\Fes,\Fes\>$, as follows from
(\ref{eq:darboux}), we obtain that $a_1+a_2\neq 0$ and
\bea
\tilde{g}&=& g-2\va\nu\Fes\\
&=& (g_1,g_2)
-\frac{2}{a_1+a_2}\left(a_1(g_1-P_1),a_2(g_2-P_2)\right)\\
&=& \frac{a_1-a_2}{a_1+a_2}\left(-g_1+\frac{2a_1}{a_1+a_2}P_1,
g_2+\frac{2a_2}{a_1+a_2}P_2\right). \eea

Now we consider $f$  as in case $(ii)$ of Theorem \ref{prop:comb2}. Again,
it suffices to prove that if $M_1$ and $M_2$ both have dimension greater
than one then any Darboux transform
$$
\tilde{g}={\cal R}_{\va,\beta}(g)\colon\, M_1\times M_2\to\R^N
$$
of a warped product $g=\Phi\,\circ(g_1\times g_2)\colon\, M_1\times
M_2\to\R^N$ of isometric immersions
$$
g_1\colon\,M_1\to\R_+^{m}\;\;\; \mbox{and}\;\;\;
g_2\colon\,M_2\to\Sf^{N-m}(1),
$$
where $\Phi\colon\,\R_+^{m}\times_{\sigma}\Sf^{N-m}(1)\to\R^N$ is an
isometry and $g_2(M_2)$ is not contained in any hypersphere of
$\Sf^{N-m}(1)$, is conformally congruent to $g$, whenever the eigenbundle
net of the  Codazzi tensor ${\cal S}_{{\va},{\beta}}=\hess \va-A_\beta$
associated to $\tilde{g}$ is the product net of $M_1\times M_2$. Set
$g_1=(h_1, \ldots, h_m)$, so that $g = (h_1,...,h_{m_1},h_m g_2)$ and
$h_m$ is the warping function of the warped product metric induced by $g$.
By Lemma~\ref{le:codtensor2}, the Codazzi tensor ${\cal
S}_{{\va},{\beta}}$  has eigenvalues
$$\lambda=A\in\R\,\,\,\mbox{and}\,\,\,\mu=A+Bh_m^{-1},\,\,\,B\neq 0.$$ Integrating $d\Fes=dg\circ{\cal S}_{{\va},{\beta}}$ and \mbox{$d{\va}=\<\Fes,dg\>$} gives
$$
\Fes=Ag+B(0,g_2)+V
$$
and
$$
\va=\frac{A}{2}\|g_1\|^2+Bh_m+\<g,V\>+c,
$$
where $V=(V_1,\ldots, V_N)\in\R^N$ and $c\in\R$.

The condition (\ref{eq:darboux})  that $g$ and $\tilde{g}$ induce
conformal metrics on $M_1\times M_2$ gives
\begin{eqnarray}\label{eq:dtransf}
\tilde{g}&=&g-2\va\nu\Fes \nonumber\\
&=& g-\frac{2h_m}{B+2Ah_m}(Ag+B(0,g_2)+V)\\
&=&\frac{1}{B+2Ah_m}
(B(h_1,\ldots,h_{m-1},-h_mg_2)-2h_mV)\nonumber
\end{eqnarray}
and \be\label{cond} \|V\|^2+2B\<g_2,V\>=2cA+\frac{BA}{2h_m}\|g_1\|^2
+\frac{B}{h_m}\<g,V\>+\frac{cB}{h_m}. \ee In particular, the preceding
equation implies that $\<g_2,V\>$ is a constant,  hence $V_j=0$ for $m\leq
j\leq N$ by the condition that $g_2(M_2)$ is not contained in any
hypersphere of $\Sf^{N-m}(1)$. Assume first that $A=0$. Then
(\ref{eq:dtransf}) reduces to
$$\tilde{g}=\bar{g}-2\frac{h_m}{B}V,$$
where $\bar{g}=(h_1,\ldots,h_{m-1},-h_mg_2)$. On the other hand,
(\ref{cond}) gives
$$\frac{h_m}{B}=\frac{1}{\|V\|^2}(\<g,V\>+c)=
\frac{1}{\|V\|^2}(\<\bar{g},V\>+c),$$ and we obtain that $\tilde{g}$ is
the composition of $g$ with the  reflection that sends $(0,g_2)$ to
$(0,-g_2)$ followed by the reflection with respect to the hyperplane
orthogonal to $V$ and the translation by the vector $-2cV/\|V\|^2$.

   Now suppose that $A\neq 0$, and set
$$
K=\frac{2}{BA}(\|V\|^2-2cA).
$$
Then (\ref{cond})  reads as \be\label{temp} \|g_1\|^2+
\frac{2}{A}\<g,V\>+\frac{2c}{A}=Kh_m. \ee Composing  $\tilde{g}$ as in
(\ref{eq:dtransf}) with the reflection that sends $(0,g_2)$ to $(0,-g_2)$
and the translation by $V/A$, we obtain $F$ that is isometric to
$\tilde{g}$ and is given by
$$
F=\frac{B}{B+2Ah_m}\left(g+\frac{V}{A}\right).
$$
On the other hand, composing $g$ with the inversion with respect to a
hypersphere with radius $r$ given by
$$
r^2=\frac{BK}{2A}=\frac{1}{A^2}(\|V\|^2-2cA)
$$
and center $-V/A$ yields
$$
G= -\frac{V}{A} + \frac{BK} {2A\|g+V/A\|^2}\left(g+\frac{V}{A}\right).
$$
Using (\ref{temp}) we obtain 
\bea \|g+V/A\|^2&=&\|g_1\|^2
+\frac{2}{A}\<g,V\>+\frac{\|V\|^2}{A^2}\\
&=& Kh_m+\frac{\|V\|^2}{A^2}-\frac{2c}{A}\\
&=& Kh_m + \frac{BK}{2A}. 
\eea Hence,
$$
G+\frac{V}{A}=\frac{B} {B+2Ah_m}\left(g+\frac{V}{A}\right)=F.\qed
$$
\vspace{2ex}

In order to complete the proof of Theorem \ref{thm:blaschke} we need the
following fact.

\begin{lemma}\po\label{curves} Let $\phi\colon\,I\subset\R\to\Q_c^2\subset \Oes$ be
a unit-speed curve with nowhere vanishing curvature $k$, where $\Oes$
denotes either $\R^2$, $\R^3$ or $\Les^3$ according as $c=0,1$ or $-1$,
respectively. Then,
\begin{itemize}
\item[$(i)$] The linear system of ODE's
\be\label{eq:s1}
\left\{\begin{array}{l} h_1'=kh_2+(A-c)h_3,\,\,\,A\in\R,\vspace{1.5ex}\\
h_2'=-kh_1,\vspace{1.5ex}\\
h_3'=h_1,
\end{array}\right.
\ee has the first integral \be\label{eq:fint}
h_1^2+h_2^2+(c-A)h_3^3=K\in\R. \ee
\item[$(ii)$] If $(h_1,h_2,h_3)$ is a solution of (\ref{eq:s1}) with initial
conditions chosen so that the constant $K$ in the right-hand-side of
(\ref{eq:fint}) vanishes, and $n$ denotes a unit normal vector to $\phi$
in $\Q_c^2$ so that $\{\phi', n\}$ is positively oriented, then
$\tilde{\phi}\colon\,I\to\Oes$ given by
$$
\tilde{\phi}=\phi-2\frac{h_3\gamma}{\<\gamma,\gamma\>},\,\,\,\,\mbox{where}\,\,\,\gamma=h_1\phi'+h_2n+ch_3\phi,
$$
is a unit-speed Ribaucour transform of $\phi$ in $\Q_c^2$.  \end{itemize}
\end{lemma}

\proof Using (\ref{eq:s1}) we obtain that the derivative of
$h_1^2+h_2^2+(c-A)h_3^3$ vanishes identically, which gives $(i)$.  Using
that $\phi''=kn-c\phi$ and $n'=-k\phi'$ we obtain
$$
\gamma'=(h_1'-kh_2+ch_3)\phi' +(h_2'+kh_1)n+c(h_3'-h_1)\phi,
$$
and hence $\gamma'=Ah_3\phi'$ by (\ref{eq:s1}). Moreover, we have
$h_3'=h_1=\<\gamma,\phi'\>$. Thus $\tilde{\phi}$ is a Ribaucour transform
of $\phi$. For $c\neq 0$, we have $\<\gamma,\phi'\>=h_3$, which implies
that $\<\tilde{\phi},\tilde{\phi}\>= \<{\phi},{\phi}\>,$ that is,
$\tilde{\phi}(I)\subset \Q_c^2$. Finally,
$$
\<\tilde{\phi}'(s),\tilde{\phi}'(s)\>
=(1-2h_3\<\gamma,\gamma\>^{-1}(Ah_3))^2=1
$$
by (\ref{eq:fint}).\qed \vspace{2ex}

Now let $f\colon\,M^n\to\R^{n+1}$ be either a cylinder
$\phi\times\R^{n-1}$ over a plane curve $\phi\colon\,I\to\R^2$, a cylinder
$C(\phi)\times\R^{n-2}$, where $C(\phi)$ is the cone over a curve
$\phi\colon\,I\to\Sf^2\subset\R^3$, or a rotation hypersurface over a
curve $\phi\colon\,I\to\R_+^2$. We prove that $f$ admits a Darboux
transform not conformally congruent to it.  Let $\tilde{\phi}\colon\,I\to
\Q_c^2$ be any Ribaucour transform of $\phi$ given by Lemma \ref{curves},
where in the case $c=-1$ we use the half-plane model of $\Hy^2$ on
$\R_+^2$. Let $\tilde{f}$ be either the cylinder
$\tilde{\phi}\times\R^{n-1}$ over $\tilde{\phi}$, the cylinder
$C(\tilde{\phi})\times\R^{n-2}$, where $C(\tilde{\phi})$ is the cone over
$\tilde{\phi}\colon\,I\to\Sf^2\subset\R^3$, or the  rotation hypersurface
over  $\tilde{\phi}\colon\,I\to\R_+^2=\Hy^2$. By Remark~\ref{re:metric}
the metrics induced by $f$ and $\tilde{f}$ are conformal, for  $\phi$ and
$\tilde{\phi}$ have the same arc-length function. Now, since $\phi$ and
$\tilde{\phi}$ are Ribaucour transforms one of each other,  there exists a
congruence of circles in $\Q_c^2$ having $\phi$ and $\tilde{\phi}$ as
envelopes. For each such circle, consider the hypersphere of $\R^{n+1}$
that intersects either $\R^2$, $\Sf^2$ or $\R_+^2=\Hy^2$ orthogonally
along it. This gives a sphere congruence  in $\R^{n+1}$ that is enveloped
by $f$ and $\tilde{f}$. We conclude that $f$ and $\tilde{f}$ are Darboux
transforms one of each other. This completes the proof of Theorem
\ref{thm:blaschke}.

\begin{remark}\po\label{re:final}
\textnormal{If we drop the assumption that $f$ and $\tilde{f}$ are not
conformally congruent  in Theorem~\ref{thm:blaschke} then we  have the
following further possibilities.
\begin{itemize}
\item[$(i)$] There exists an inversion $I$ in
$\R^{n+1}$ such that $L'(\tilde{f})=I(L(f))$, where $L$ and $L'$ are
compositions of a homothety and a translation.
\item[$(ii)$] The immersions $f$ and $\tilde{f}$ are as in either case of Theorem \ref{prop:comb2}
with both factors of dimension greater than one.
\end{itemize} }
\end{remark}

\vspace{.5in} {\renewcommand{\baselinestretch}{1}

\hspace*{-20ex}\begin{tabbing} \indent\= IMPA -- Estrada Dona Castorina, 110
\indent\indent\=  Universidade Federal de S\~{a}o Carlos\\
\> 22460-320 -- Rio de Janeiro -- Brazil  \>
13565-905 -- S\~{a}o Carlos -- Brazil \\
\> E-mail: marcos@impa.br \> E-mail: tojeiro@dm.ufscar.br
\end{tabbing}}
\end{document}